\documentclass[a4paper,11pt]{amsart}
\usepackage{amsmath,amsthm,amssymb}
\usepackage{comment}
\newcommand{\mychoice}[2]{#1
}
\mychoice{
\newcommand{\plabel}[1]{ \label{#1}}
\newcommand{\gbibitem}[1]{ \bibitem{#1}}
\newcommand{\snewpage}{}

\specialcomment{commentx}{}{}\excludecomment{commentx}
}{
\usepackage{xcolor}
\newcommand{\plabel}[1]{ \label{#1}\rlap{\smash{${}^{^{[#1]}}$}}}
\newcommand{\gbibitem}[1]{ \bibitem{#1}\rlap{\smash{${}^{^{[#1]}}$}}}
\newcommand{\snewpage}{\newpage}

\newenvironment{commentx}{\color{magenta} }{\color{black} }
}
\DeclareMathOperator{\supp}{supp}
\DeclareMathOperator{\dec}{dec}
\DeclareMathOperator{\proj}{proj}
\DeclareMathOperator{\conv}{conv}
\DeclareMathOperator{\pgen}{pgen}
\DeclareMathOperator{\rgen}{rgen}
\DeclareMathOperator{\relab}{rlb}
\DeclareMathOperator{\ext}{ext}
\DeclareMathOperator{\cmon}{cmon}
\DeclareMathOperator{\card}{card}
\DeclareMathOperator{\reg}{reg}
\DeclareMathOperator{\topp}{top}
\DeclareMathOperator{\cand}{cand}
\DeclareMathOperator{\alt}{alt}
\DeclareMathOperator{\forw}{forw}

\theoremstyle{definition}
\newtheorem{point}{}[section]

\newtheorem{remark}[point]{Remark}
\theoremstyle{plain}

\newtheorem{theorem}[point]{Theorem}
\newcommand{\marginextend}[1]{ \addtolength{\oddsidemargin}{-#1}  \addtolength{\evensidemargin}{-#1}
  \addtolength{\textwidth}{#1}\addtolength{\textwidth}{#1}}
\newcommand{\updownextend}[1]{ \addtolength{\topmargin}{-#1}  \addtolength{\textheight}{#1}
\addtolength{\textheight}{#1}}
\marginextend{1cm}
\updownextend{0cm}
\allowdisplaybreaks[4]
\newcommand{\bem}{\begin{bmatrix}}
\newcommand{\eem}{\end{bmatrix}}

\title{Alternating patterns in commutator monomials}
\author{Gyula Lakos}
\address{Department of Geometry, Institute of Mathematics, E\"otv\"os University,
P\'azm\'any P\'eter s.~1/C,  Budapest, H--1117, Hungary}
\email{lakos@cs.elte.hu}
\begin{document}
\begin{abstract}
Considering commutator monomials of the non-commutative
associative variables $X_1,\ldots,X_n$; we determine the maximal possible number of
alternating associative monomials in their noncommutative polynomial expansions.
This is achieved by replacing generating functions with polytope sequences, which turn out to be finitely generated in some sense.
\end{abstract}
\maketitle
\snewpage
\section{Introduction and statement of main results}
We call a permutation $\sigma$ of $\{1,\ldots,n\}$ alternating,
if $\sigma(k)$ is not between $\sigma(k-1)$ and $\sigma(k+1)$ for any  $k\in\{2,\ldots,n-1\}$.
This terminology is somewhat different from the one of  Flajolet, Sedgewick \cite{FS} or Stanley \cite{Sta},
whose ``forward'' alternating permutations are required to do not start with ascents%
%, and reverse alternating permutations do not start with descents
.

The following question was originally motivated by investigations concerning
the convergence of the Magnus expansion in Banach--Lie algebras, cf.~\cite{L3}.
Assume that $X_1,\ldots,X_n$ are noncommutative associative variables.
Let us consider the commutator monomials $M$ built up from $X_1,\ldots,X_n$ by $n-1$ many ``Lie'' commutators
(i.~e.~commutators $\boldsymbol[Y,Z\boldsymbol]=YZ-ZY$), with all variables used; e.~g.~$M=[[X_5,[X_1,X_4]],[X_3,X_4]]$ for $n=5$.
We ask what is the maximal number $r_n$
of  associative monomials $X_{\sigma(1)}\cdot\ldots\cdot X_{\sigma(n)}$ in the expanded $M$ such that $\sigma$ is alternating.
Trivially, the answer is the same if we ask the question about  ``Jacobi'' commutators
(i.~e.~commutators $\boldsymbol\{Y,Z\boldsymbol\}=YZ+ZY$).
(The reader is free to reformulate the question in terms of binary trees, purely in combinatorial terms.)
Although, ultimately, the answer turns out to be not particularly informative for the purposes the of Magnus expansion,
the question is a strangely engaging one.

It is not hard too see that the expected number of
alternating associative monomials in the commutators monomials above
(averaged over permutations of the variables $X_1,\ldots,X_n$) is asymptotically $\sim (4/\pi)^{n+1}$, cf. Appendix \ref{sec:asy}.
Thus, asymptotically, we expect at least that much for $r_n$.
The exact answer is
\begin{theorem}
\plabel{th:rn}
The number $r_n$ is given by the table
\[
\begin{array}{c||c|c|c|c|c|c|c|c|c|}
n&1&2&3&4&5&6+4k&7+4k&8+4k&9+4k\\
\hline
r_n&1&2&4&6&8&14\cdot6^k&24\cdot6^k&36\cdot6^k&56\cdot6^k
\end{array}
\]
where $k\in\mathbb N$.
\end{theorem}

In particular, this is asymptotically $c_{n\,\mathrm{mod}\,4}\cdot (6^{1/4})^n$;
which is not that much greater than the averaged value.
(Being $4/\pi= 1.27323\ldots$ and $6^{1/4}=1.56508\ldots$ )

Although it is not our primary interest, we may also ask
about the maximal possible number $r_n^{\forw}$ of monomials of ``forward'' alternating type
in commutator monomials.
\begin{theorem}
\plabel{th:rpn}
The number $r_n^{\forw}$ is given by the table
\[
\begin{array}{c||c|c|c|c|c|c|c|}
n&1&2&3&4+4k&5+4k&6+4k&7+4k\\
\hline
r_n^{\forw}&1&1&2&3\cdot6^k&6\cdot6^k&7\cdot6^k&14\cdot6^k
\end{array}
\]
where $k\in\mathbb N$.
\end{theorem}
Similar comment applies about the asymptotic behaviour, except  the averaged expected number is halved (for $n\geq2$).

As the exact answers above are rather simple (with only $n\in\{1,2,5\}$ falling out of the general pattern
in Theorem \ref{th:rn}, and only $n\in\{2,3\}$ being irregular in Theorem \ref{th:rpn}),
one might hope that the theorems above can be proven by induction on $n$.
This is sort of true, but, in practice, one gets quickly frustrated by the large number
of combinatorial possibilities one faces even in low orders.
A way to deal with the combinatorial multitude of the problem is, however, provided by convex geometry / linear programming.
The objective of this paper is to explain this method.
The basic is idea is to replace generating functions with polytope sequences, which turn out to be finitely generated in some sense.

In order to do so, we extend our terminology.
Firstly, we will relax the requirement that in the commutator monomials we use
exactly the variables $X_1,\ldots,X_n$ with multiplicity $1$.
Instead, we allow to use arbitrary variables $X_i$ with $i\in\mathbb N^\star\equiv\mathbb N\setminus\{0\}$.
In general, a function $\sigma:\{1,\ldots,n\}\rightarrow \mathbb N^\star$ will be called alternating
if, for any $k\in\{1,\ldots,n-1\}$, it holds that $\sigma(k)\neq\sigma(k+1)$,
and, for any $k\in\{2,\ldots,n-1\}$, it holds that $\sigma(k)$ is not between $\sigma(k-1)$ and $\sigma(k+1)$.
If $P$ is a homogeneous polynomial of order $n$ of some noncommutative variables $X_i$ ($i\in\mathbb N^\star$),
\begin{equation}
P=\sum_{\sigma:\{1,\ldots,n\}\rightarrow \mathbb N^\star}p_\sigma X_{\sigma(1)}\ldots X_{\sigma(n)},
\plabel{eq:P}
\end{equation}
then we let
\[\|P\|_{\alt}:=\sum_{\substack{\sigma:\{1,\ldots,n\}\rightarrow \mathbb N^\star \\ \text{is alternating} }}|p_\sigma |.\]
We can ask what is the maximal number of $\|M\|_{\alt}$, where
$M$ is a commutator monomial of order $n$ of some variables $X_i$ ($i\in\mathbb N^\star$).
The answer is $r_n$ again.
Indeed, if $M$ was a commutator monomial of the variables $X_{i_1},\ldots,X_{i_n}$ with
$i_1\leq\ldots\leq i_n$, then we can consider a commutator monomial $M^{\mathrm c}$ which is like $M$
but the variables are replaced by $X_1,\ldots,X_n$.
Then $\|M\|_{\alt}\leq \|M^{\mathrm c}\|_{\alt}$, and $ \|M^{\mathrm c}\|_{\alt}$ is just the number
of alternating associative monomials in $M^{\mathrm c}$.
This argument applies for Lie commutators and Jacobi commutators alike.
We may call this process as ``combing out'' the multiplicities.
Note that this process is not deterministic: the result may vary because the variables
with multiplicities can be replaced in multiple ways.
%Moreover this ``combing out'' can conveniently made wider by replacing $X_{i_1},\ldots,X_{i_n}$
%by  $X_{j_1},\ldots,X_{j_n}$ where $j_1<\ldots<j_n$.

Let us fix here that we use only Jacobi commutators in the future,
and we consider only polynomials \eqref{eq:P} with $p_\sigma\geq0$.
Then we do not have to keep the alternating associative monomials in their individual form.
For example, the associative monomial $X_2X_5X_2X_8X_1$ can be replaced
by the symbol $X_{2\nearrow\searrow1}$, keeping only the first and last indices, and the information that they are followed or preceded by ascents or descents.
Similarly, the associative monomial $X_2X_5X_6X_8X_1$ can be replaced by $0$ as it is not an alternating monomial.
In that way we obtain the linear factor space $\mathcal X$ generated linearly by the terms
$\tilde X_i$, $X_{i\nearrow\nearrow j}$,  $X_{i\searrow\searrow j}$, $X_{i\nearrow\searrow j}$, $X_{i\searrow\nearrow j}$,
where $i,j\in\mathbb N^\star$.
Note that the Lie and Jacobi commutators  act naturally on $\mathcal X$.
However, in the Jacobi commutator case (with nonnegative coefficients),
$\|\cdot\|_{\alt}$ is compatible with the factorization process.
Ultimately, on $\mathcal X^{\geq0}$, we have a natural action of the Jacobi commutator, and
an induced norm $\|\cdot\|_{\alt}$ (old notation is kept) with respect to additivity and nonnegative scaling.
Written explicitly, if
\begin{multline}
\tilde P=\sum_{i\in\mathbb N^\star}p_i\tilde X_i
+\sum_{i,j\in\mathbb N^\star}\bigl(
p_{i\nearrow\nearrow j}X_{i\nearrow\nearrow j}
+\\+p_{i\searrow\searrow j}X_{i\searrow\searrow j}
+p_{i\nearrow\searrow j}X_{i\nearrow\searrow j}
+p_{i\searrow\nearrow j}X_{i\searrow\nearrow j}
\bigr)
\plabel{eq:PT}
\end{multline}
with nonnegative coefficients (i.~e.~it is in $\mathcal X^{\geq0}$), then
\[\|\tilde P\|_{\alt}=
\sum_{i\in\mathbb N^\star}p_i
+\sum_{i,j\in\mathbb N^\star}\left(
p_{i\nearrow\nearrow j}
+p_{i\searrow\searrow j}
+p_{i\nearrow\searrow j}
+p_{i\searrow\nearrow j}
\right).
\]
(Algebraically, with respect to the commutators,
  $\tilde X_i=\tilde X_{i\nearrow\searrow i}+\tilde X_{i\searrow\nearrow i}$,
where $\tilde{}$ is lost after being commutatored.)
If $\tilde P$ is as in \eqref{eq:PT}, then we let
\[\supp
\tilde P=\bigcup_{\substack{i\in\mathbb N^\star\\ p_i\neq0}}\{i\}
\cup\bigcup_{\substack{i,j\in\mathbb N^\star\\
p_{i\nearrow\nearrow j}\neq0 \text{ or }
p_{i\searrow\searrow j}\neq0 \text{ or }
p_{i\nearrow\searrow j}\neq0 \text{ or }
p_{i\searrow\nearrow j}\neq0
 }}\{i,j\};
\]
i.~e.~the set of indices which occur in the expansion of $\tilde P$.
If $f:\supp \tilde P\subset I\subset \mathbb N^\star\rightarrow\mathbb N^\star$ is a (not necessarily strictly) monotone increasing function, then we let
\begin{multline}
\pi_f\tilde P=\sum_{i\in\mathbb N^\star}p_i\tilde X_{f(i)}
+\sum_{i,j\in\mathbb N^\star}\bigl(
p_{i\nearrow\nearrow j}X_{f(i)\nearrow\nearrow f(j)}
+\\+p_{i\searrow\searrow j}X_{f(i)\searrow\searrow f(j)}
+p_{i\nearrow\searrow j}X_{f(i)\nearrow\searrow f(j)}
+p_{i\searrow\nearrow j}X_{f(i)\searrow\nearrow f(j)}
\bigr).
\plabel{eq:PTf}
\end{multline}
Note  that $\|\tilde P\|_{\alt}=\|\pi_f\tilde P\|_{\alt}$.
If $f$ is injective on the support of $\tilde P$, then $\pi_f\tilde P$ is called a relabeling of $P$.
If $f$ is defined on the support of $\tilde P$, sending the support
to a set of shape $\{1,\ldots,k\}$, then $\pi_f\tilde P$ is called a regular projection of $P$.
Note that any $\tilde P$ has typically multiple but only finitely many regular projections.
We have a regular relabeling when the support $\{i_1,\ldots,i_k\}$ (with $i_s<i_{s+1}$) is sent to $\{1,\ldots,k\}$.
The regular relabeling $\reg \tilde P$ is unique.
In general, any projection $\pi_f$ can be imagined as a succession of a regular projection and a relabeling.
If $\{i_1,\ldots,i_k\}$ (with $i_s<i_{s+1}$) is the support of $P$, and $f$ sends $i_s$ to $j_s$, then we may use the notation
\[\pi_f\tilde P\equiv\pi^{[i_1,\ldots,i_k]}_{[j_1,\ldots,j_k]}\tilde P.\]

Assume that $\mathcal S\subset \mathcal X^{\geq0}$.
We let (the projected set) $\proj\mathcal S$ be the set of all $\pi_f x$,
where $x\in\mathcal S$ and $f$ is defined on the support of $x$.
We let (the decreased set) $\dec\mathcal S$ be the set of those $y\in \mathcal X^{\geq0}$ for which there exits an element $x\in \mathcal S$
such that $0\leq y\leq x$ (monomial-wise).
Furthermore, let  $\conv\,\mathcal S$ be the set of convex combinations from $\mathcal S$.
Then, we define the projective generatum $\pgen \mathcal S$ of $\mathcal S$ as $\conv\dec\proj\mathcal S$.
It is easy to see that $\pgen \mathcal S$ will be closed for the operations $\conv,$ $\dec,$ $\proj$.
Now, $\proj$ can be decomposed into two steps:
In general, we let (the regularly projected set) $\proj_0\mathcal S$ be the set of all $\pi_f x$, where $\pi_f x$ is a regular projection,
and we  let (the relabeled set) $\relab\mathcal S$ be the set of all $\pi_f x$, where $\pi_f x$ is a relabeling.
Thus $\proj\mathcal S=\relab\proj_0\mathcal S$.
We define the relatively linear generatum $\rgen \mathcal S$ of $\mathcal S$ as $\conv\dec\relab\mathcal S$.
Then
\[\pgen\mathcal S=\rgen \proj_0\mathcal S.\]
All these operations leave $\|S\|_{\alt}=\sup_{x\in\mathcal S}\|x\|_{\alt}$ invariant.
Also, the operations above are compatible to the Jacobian commutators in the following sense:
the Jacobian commutator of convex combinations of is a convex combination of Jacobian commutators;
the Jacobian commutator of  decreased terms of is a decreased term of the Jacobian commutator of the original terms;
the Jacobian commutator of projections is a decreased term of a projection of a Jacobian commutator of some relabelings of the original terms
with disjoint supports (this is a case for combing out some multiplicities).

Let $\mathcal B_n^{\cmon}$ be the set all (Jacobian) commutator monomials
built up from the starting elements $\tilde X_1,\ldots,\tilde X_n$ using $n-1$ many commutators (with all starting elements used).
Then $\mathcal B_n^{\cmon}$ contains finitely many elements with $r_n=\|\mathcal B_n^{\cmon}\|_{\alt}$.
We define
\[\mathcal B_n=\pgen \mathcal B_n^{\cmon}.\]
Then we have $r_n=\|\mathcal B_n\|_{\alt}$.
Also note that $\mathcal B_n=\rgen \left(\proj_0 \mathcal B_n^{\cmon}\right)$,
where $\proj_0 \mathcal B_n^{\cmon}$ is still finite.
The problem is that the generating systems above for $\mathcal B_n$ are quite large.
Yet, we may hope that $\mathcal B_n$ allows simpler generating systems.
And this is, indeed, the case.
In order to formulate the corresponding result, we define the normalized sets
\[\mathcal R_n=\frac1{(6^{1/4})^n}\mathcal B_n.\]
Now we can state
\begin{theorem}
\plabel{th:Rn}
\[\mathcal R_2\subsetneqq\mathcal R_6\subsetneqq\mathcal R_{10}=\mathcal R_{14}=\ldots\quad\text{(stabilizes);}\]
\[\mathcal R_3\subsetneqq\mathcal R_7\subsetneqq\mathcal R_{11}\subsetneqq\mathcal R_{15}=\mathcal R_{19}=\ldots\quad\text{(stabilizes)};\]
\[\mathcal R_4\subsetneqq\mathcal R_8\subsetneqq\mathcal R_{12}=\mathcal R_{16}=\ldots\quad\text{(stabilizes)};\]
\[\mathcal R_5\subsetneqq\mathcal R_9=\mathcal R_{13}=\ldots\quad\text{(stabilizes).}\]
\end{theorem}
Using Theorem \ref{th:Rn},  Theorems \ref{th:rn}  and \ref{th:rpn} are ``easy'' to check:
It is sufficient to verify $n\in\{1,\ldots,12,15\}$; and the rest follows by the stabilization property.
This, however, still involves checking several cases, and Theorem \ref{th:Rn} is expectedly even harder.
What makes the situation tractable is that there are generating systems of tractable size for $\mathcal R_i$.
The content of the next section is the
description of these generating systems and the  outline of the proofs of the above mentioned theorems.

\snewpage
\section{Commutator polytopes and their generating sets}
Before describing the various $\mathcal B_n$ or $\mathcal R_n$  by small generating sets, let us consider
what kind of nice generating sets we might hope for.
Here we use some convex geometry; see Matou\v{s}ek, G\"{a}rtner \cite{MG} for more about that.

Firstly, let us consider the relatively linear ($\rgen$-) generating sets for $\mathcal B_n$.
We know that $\mathcal B_n$ is relatively linearly generated by $\proj_0\mathcal B_n^{\cmon}$, which is finite set,
and regular in the sense that the supports are regular (i.~e.~of shape $\{1,\ldots,k\}$).
We will see that there is a single minimal regular relatively linear generating set for $\mathcal B_n$,
which is then necessarily is a subset of $\proj_0\mathcal B_n^{\cmon}$.
In general, for a finite set $\mathcal I\subset \mathbb N^\star$, let
\[\mathcal S|_{\mathcal I}=\{x\in\mathcal S\,:\,\supp x\subset I\}.\]
For later use, we also introduce the notation
\[\mathcal S|_{\reg}=\{x\in\mathcal S\,:\,x=\reg x\}.\]
Then it is easy to see that
\[\mathcal B_n|_{\mathcal I}=\conv\dec\proj_{\mathcal I}\mathcal B_n^{\cmon},\]
where, in general, $\proj_{\mathcal I}\mathcal S$ is the set of projections where the support of the elements of
$\mathcal S$ is sent to a subset of $\mathcal I$.
Here, $\mathcal B_n|_{\mathcal I}\subset \mathcal X^{\geq0}|_{\mathcal I}$, the latter one being a finite dimensional
positive space.
Now $\conv^+=\conv\dec$ is a kind of convex closure on $\mathcal X^{\geq0}|_{\mathcal I}$, for which
there is a minimal generating set $\ext^+ \proj_{\mathcal I}\mathcal B_n^{\cmon}$
such that
\[\conv\dec\proj_{\mathcal I}\mathcal B_n^{\cmon}=\conv\dec\ext^+\proj_{\mathcal I}\mathcal B_n^{\cmon}.\]
In fact, in the dual viewpoint $\ext^+$ is set points which are singular extrema for strictly positive linear functionals on
$\mathcal X^{\geq0}|_{\mathcal I}$.
Here the finitely generatedness was crucial.
In any case,
\[\ext^+\mathcal B_n|_{\mathcal I}=\ext^+ \proj_{\mathcal I}\mathcal B_n^{\cmon}.\]
(And note that $\ext^+$ is always a selection of elements.)
Up to this point we have used only convex geometry and finitely generatedness.
What is less trivial is that for $\mathcal I \subset \mathcal I'$
\[\ext^+ \proj_{\mathcal I}\mathcal B_n^{\cmon}\subset \ext^+ \proj_{\mathcal I'}\mathcal B_n^{\cmon}\]
holds, i.~e.~there is no loss of extremal points. (Indeed, for
$x\in(\ext^+ \proj_{\mathcal I}\mathcal B_n^{\cmon})\setminus(\ext^+ \proj_{\mathcal I'}\mathcal B_n^{\cmon})$,
due to the support,
the element would be strictly majorized by an element $z$ of  $\mathcal B_n|_{\mathcal I'}$, but then a
projection of $z$ to $\mathcal I$ would strictly majorize $x$, which is a contradiction.)
Due to relabelings, it is clear that in the structure of $\mathcal B_n|_{\mathcal I}$
only the cardinality of $\mathcal I$ is important.
Thus it is sufficient to consider the cases $\mathcal I=\{1,\ldots,k\}$.
We have a sequence of sets
\begin{equation}
\ext^+ \proj_{\{1\}}\mathcal B_n^{\cmon}\subset\ldots\subset
\ext^+ \proj_{\{1,\ldots,k\}}\mathcal B_n^{\cmon}\subset\ldots
\plabel{eq:hier}
\end{equation}
As extremal points will be invariant for relabelings,
these sets can be characterized by the regular restrictions
\begin{equation}
\left.\left(\ext^+ \proj_{\{1\}}\mathcal B_n^{\cmon}\right)\right|_{\reg}\subset\ldots\subset
\left.\left(\ext^+ \proj_{\{1,\ldots,m\}}\mathcal B_n^{\cmon}\right)\right|_{\reg}
\plabel{eq:acon}
\end{equation}
(and the knowledge of the full support sets $\{1,\ldots, k\}$).
This latter sequence will stabilize due finitely generatedness by $\left.\left(\ext^+  \proj_{\{1,\ldots,n\}}\mathcal B_n^{\cmon}\right)\right|_{\reg}$
(as further extensions introduce only relabelings).
Here $m$ is the smallest value for which the sequence stabilizes, and we know that $m\leq n$.
Every generating set of  $\mathcal B_n$ with respect to $\conv^+\equiv \conv \dec$
must contain all relabelings of $ \left.\left(\ext^+ \proj_{\{1,\ldots,m\}}\mathcal B_n^{\cmon}\right)\right|_{\reg}$;
and, conversely, the relabelings already form a generating set for  $\mathcal B_n$.
Moreover, every regular generating set of
 $\mathcal B_n$ with respect to $\rgen\equiv \conv \dec\relab$
must contain $\left.\left(\ext^+ \proj_{\{1,\ldots,m\}}\mathcal B_n^{\cmon}\right)\right|_{\reg}$, which already forms a regular relatively linear
generating set.
We emphasize that the hierarchy of extremal sets  depend only on $\mathcal B_n$ as
\begin{equation}
\ext^+\mathcal B_n|_{\{1,\ldots,k\}}=\ext^+ \proj_{\{1,\ldots,k\}}\mathcal B_n^{\cmon};
\plabel{eq:elon}
\end{equation}
and for the regular extremal points
\[
\left.\left( \ext^+\mathcal B_n|_{\{1,\ldots,k\}}\right)\right|_{\reg}
=
\left.\left( \ext^+ \proj_{\{1,\ldots,k\}}\mathcal B_n^{\cmon}\right)\right|_{\reg}
.\]
In general, $\mathcal B_n^{\cmon}$ can be replaced by any finite projective generating set $\mathcal G_n$ of $\mathcal B_n$.
In that case, we know  $m\leq \max \card \supp \mathcal G_n$ (maximal cardinality of support); which information
might be of advantage, as we might have to compute the extremal sets in fewer cases.

Next we consider projective ($\pgen$-) generating sets for $\mathcal B_n$.
Here $\mathcal B_n=\pgen \mathcal S$ if and only if $\mathcal S\subset\mathcal B_n$
and $\ext^+ \proj_{\{1,\ldots,m\}}\mathcal B_n^{\cmon}\subset\pgen \mathcal S$.
Now, $\mathcal S=\mathcal B_n^{\cmon}$ fits this role, but according to Theorem \ref{th:Rn}
we cannot hope for a single minimal projective generating set for $\mathcal B_n$, and there are ones possibly
disjoint from $\mathcal B_n^{\cmon}$.
Let
\begin{equation}
\topp_\pi \reg \ext^+ \proj_{\{1,\ldots,m\}}\mathcal B_n^{\cmon}
\plabel{eq:top}
\end{equation}
be denote the set those elements $x\in \reg \ext^+ \proj_{\{1,\ldots,m\}}\mathcal B_n^{\cmon}$ for which there
are no $y\in \reg \ext^+ \proj_{\{1,\ldots,m\}}\mathcal B_n^{\cmon}$ such that $x=\pi_fy$ with some $f$.
Due to finiteness and the natural filtration by the support, we see that
\eqref{eq:top} yields a projective generating system for $\mathcal B_n$.
Again, we note that this generating system depends only on $\mathcal B_n$, canonically.

We define the cleaning operation $T:\mathcal X^{\geq0}\rightarrow\mathcal X^{\geq0}$
which simply replaces the terms $X_{i\nearrow\searrow i}$, $X_{i\searrow\nearrow i}$,
$X_{i\nearrow\nearrow i}$, $X_{i\searrow\searrow i}$ by $0$.
We say that $x\in\mathcal X^{\geq0}$ is clean if $Tx=x$.
In particular, the elements of $\mathcal B_n^{\cmon}$ are clean:
Assume that $x\in\ext^+ \proj_{\{1,\ldots,m\}}\mathcal B_n^{\cmon}$ is not clean.
We know that there is an element $y\in \mathcal B_n^{\cmon}$ such that $x=\pi_fy$ with some $f$.
Then
\begin{equation}
y\leq c_1x_1+\ldots+c_sx_s,
\plabel{eq:up}
\end{equation}
where $x_i\in\ext^+ \proj_{\{1,\ldots,n\}}\mathcal B_n^{\cmon}$ and $c_i>0$ with $c_1+\ldots+c_s=1$.
Applying $\pi_f$, we see that
\begin{equation}
x\leq c_1\pi_fx_1+\ldots+c_s\pi_fx_s.
\plabel{eq:down}
\end{equation}
Then by the extremality of $x$ we see that $x=\pi_f x_i$ and there is equality in \eqref{eq:down}.
As projection leaves $\|\cdot\|_{\alt}$ invariant, this forces equality in \eqref{eq:up}.
As $y$ is clean, this forces the $x_i$ to be clean.
As $x=\pi_f x_i$, we see that $x$ is $\pi$-topped by a clean $x_i$.
In general, this argument shows that \eqref{eq:top} contains only clean elements.

This observation helps a little bit to recover \eqref{eq:top}, as it is sufficient to consider clean elements.
We may consider the sets
\begin{equation}
\ext^+T\mathcal B_n|_{\mathcal I}=\ext^+ T\proj_{\mathcal I}\mathcal B_n^{\cmon};
\plabel{eq:elon2}
\end{equation}
(restriction to $\mathcal I$ and cleaning commute), which leads to the sequence
\begin{equation}
\ext^+ T\proj_{\{1\}}\mathcal B_n^{\cmon},\ldots,
\ext^+ T\proj_{\{1,\ldots,k\}}\mathcal B_n^{\cmon},\ldots
.
\plabel{eq:hier2}
\end{equation}
(The first set is empty, of course.)
This is simpler than \eqref{eq:hier} in the sense that
\[\ext^+ T\proj_{\mathcal I}\mathcal B_n^{\cmon}=\ext^+T\ext^+ \proj_{\mathcal I}\mathcal B_n^{\cmon},\]
in general; thus the structure of the individual sets is reduced.
On the other hand, \eqref{eq:elon2} is not necessarily monotone in $\mathcal I$, and
it is not necessarily closed for relabelings inside $\mathcal I$.
In particular, \eqref{eq:hier2} cannot be characterized simply by the regular restriction.

For this reason, in general, let
\[\ext^{+0}\mathcal S=\bigcup_{\mathcal J}\ext^+ \mathcal S|_{\mathcal J}.\]
Then
\begin{equation}
\ext^{+0}T\mathcal B_n|_{\mathcal I}=\ext^{+0}  T\proj_{\mathcal I}\mathcal B_n^{\cmon}.
\plabel{eq:elon3}
\end{equation}
This is clearly monotone in $\mathcal I$. The sequence
\begin{equation}
\ext^{+0} T\proj_{\{1\}}\mathcal B_n^{\cmon}\subset\ldots\subset
\ext^{+0} T\proj_{\{1,\ldots,k\}}\mathcal B_n^{\cmon}\subset\ldots
\plabel{eq:hier3}
\end{equation}
is more similar to \eqref{eq:hier}; its members are closed for relabelings (in the indicated support set),
thus they can be characterized by the regular restriction.
In fact,
\begin{multline}
\left.\left(\ext^{+0}T\mathcal B_n|_{\{1,\ldots,k\}} \right)\right|_{\reg}=
\left.\left(\ext^{+0} T\proj_{\{1,\ldots,k\}}\mathcal B_n^{\cmon}\right)\right|_{\reg}=
\\
=\left.\left(\ext^+ T\proj_{\{1\}}\mathcal B_n^{\cmon}\right)\right|_{\reg}
\cup\ldots\cup
\left.\left(\ext^+ T\proj_{\{1,\ldots,k\}}\mathcal B_n^{\cmon}\right)\right|_{\reg}
.\end{multline}
Here the stabilizing threshold will be the same $k=m$, the maximal cardinality for the supports in \eqref{eq:top}.
Then
%To keep the long story short,
\begin{equation}
%\begin{multline}
\topp_\pi \left.\left( \ext^+ \proj_{\{1,\ldots,m\}}\mathcal B_n^{\cmon}\right)\right|_{\reg}
%=
%\topp_\pi \left.\left( \ext^+ T\proj_{\{1,\ldots,m\}}\mathcal B_n^{\cmon}\right)\right|_{\reg}
%\\
=\topp_{T\pi}
\left.\left(\ext^{+0} T\proj_{\{1,\ldots,m\}}\mathcal B_n^{\cmon}\right)\right|_{\reg};
\plabel{eq:titop}
%\end{multline}
\end{equation}
where the latter set contains those elements $x\in \left.\left(\ext^{+0} T\proj_{\{1,\ldots,m\}}\mathcal B_n^{\cmon}\right)\right|_{\reg}$
such that there is no element $y\in \left.\left(\ext^{+0} T\proj_{\{1,\ldots,m\}}\mathcal B_n^{\cmon}\right)\right|_{\reg}$
with the property $x=T\pi_fy$ with some $f$.
As we have seen, via \eqref{eq:elon2}, the extremal points in \eqref{eq:elon3} are all cleaned versions of \eqref{eq:elon}
except possibly not all of those.
This looks like a small gain, but in a computationally heavy environment it counts a bit.
Again, the projective generating set $\mathcal B_n^{\cmon}$ of $\mathcal B_n$ can be replaced by any other finite $\mathcal G_n$.
Altogether, we find:

The theoretically ideal way to describe $\mathcal B_n$ is by the ``full canonical generating system''
\[\left.\left( \ext^+ \mathcal B_n\right)\right|_{\reg}=  \left.\left( \ext^+ \proj_{\{1,\ldots,m\}}\mathcal G_n\right)\right|_{\reg}.\]
A cruder but more practical way is by the ``reduced canonical generating system''
\[\left.\left(\ext^{+0} T\mathcal B_n\right)\right|_{\reg}=  \left.\left( \ext^{+0}T \proj_{\{1,\ldots,m\}}\mathcal G_n\right)\right|_{\reg}.\]
Finally, a concise way to describe $\mathcal B_n$ is by the ``top canonical generating system''
\[\topp_\pi\left.\left( \ext^+ \mathcal B_n\right)\right|_{\reg}= \topp_{T\pi}\left.\left(\ext^{+0} T\mathcal B_n\right)\right|_{\reg}. \]
These are all natural projective generating sets for $\mathcal B_n$ not depending on choice of $\mathcal G_n$.

Now, $\mathcal R_n$ is the same as $\mathcal B_n$ but multiplied by $6^{-n/4}$.
%It is clear that
$\mathcal R_1$ is projectively generated by
\[R^{[1]}_1=\frac1{6^{1/4}}X_1. \]
For the sake of the next statement, we introduce the notation
\[
\lfloor A,B\rfloor=\frac1{6^{1/4}}\boldsymbol\{A,B\boldsymbol\}
.\]
\snewpage
\begin{theorem}
\plabel{th:gen}
(a) The top canonical generating system for $\mathcal R_2$ is
given by the single element
\[ R_{{{1}}}^{[2]}= \pi_{{{[1,2]}}}^{[1,2]}\lfloor   \pi_{{{[1]}}}^{[
1]}   R_{{{1}}}^{[1]}, \pi_{{{[2]}}}^{[1]}   R_{{{1}}}^{[1]}\rfloor
.\]

(b) The top canonical generating system for $\mathcal R_3$ is
given by the single element
\[ R_{{{1}}}^{[3]}= \pi_{{{[1,2,3]}}}^{[1,2,3]}\lfloor   \pi_{{{[2]}}
}^{[1]}   R_{{{1}}}^{[1]}, \pi_{{{[1,3]}}}^{[1,2]}   R_{{{1}}}^{[
2]}\rfloor
.\]

(c) The top canonical generating system for $\mathcal R_4$ is
given by $R_{{{1}}}^{[4]},R_{{{2}}}^{[4]} $; where:
\[ R_{{{1}}}^{[4]}= \pi_{{{[1,2,3]}}}^{[2,3,4]}\lfloor   \pi_{{{[2]}}
}^{[1]}   R_{{{1}}}^{[1]}, \pi_{{{[1,3,4]}}}^{[1,2,3]}   R_{{{1}}
}^{[3]}\rfloor\]
\[ R_{{{2}}}^{[4]}= \pi_{{{[1,2,3]}}}^{[1,2,3]}\lfloor   \pi_{{{[3]}}
}^{[1]}   R_{{{1}}}^{[1]}, \pi_{{{[1,2,4]}}}^{[1,2,3]}   R_{{{1}}
}^{[3]}\rfloor\]

(d) The top canonical generating system for $\mathcal R_5$ is
given by $R_{{{1}}}^{[5]},\ldots, R_{{{8}}}^{[5]}$; where:
\[ R_{{{1}}}^{[5]}= \pi_{{{[1,2,2]}}}^{[1,2,3]}\lfloor   \pi_{{{[1]}}
}^{[1]}   R_{{{1}}}^{[1]}, \pi_{{{[2,3,4]}}}^{[1,2,3]}   R_{{{1}}
}^{[4]}\rfloor\]
\[ R_{{{2}}}^{[5]}= \pi_{{{[1,1,2]}}}^{[2,3,4]}\lfloor   \pi_{{{[4]}}
}^{[1]}   R_{{{1}}}^{[1]}, \pi_{{{[1,2,3]}}}^{[1,2,3]}   R_{{{2}}
}^{[4]}\rfloor\]
\[ R_{{{3}}}^{[5]}= \pi_{{{[1,2,3,4]}}}^{[1,2,3,4]}\lfloor   \pi_{{{[
2]}}}^{[1]}   R_{{{1}}}^{[1]}, \pi_{{{[1,3,4]}}}^{[1,2,3]}   R_{{
{2}}}^{[4]}\rfloor\]
\[ R_{{{4}}}^{[5]}= \pi_{{{[1,2,3,4]}}}^{[1,2,3,4]}\lfloor   \pi_{{{[
3]}}}^{[1]}   R_{{{1}}}^{[1]}, \pi_{{{[1,2,4]}}}^{[1,2,3]}   R_{{
{1}}}^{[4]}\rfloor\]
\[ R_{{{5}}}^{[5]}= \pi_{{{[1,2,3,4,5]}}}^{[1,2,3,4,5]}\lfloor
\pi_{{{[1,4]}}}^{[1,2]}   R_{{{1}}}^{[2]}, \pi_{{{[2,3,5]}}}^{[1,2,3]}
   R_{{{1}}}^{[3]}\rfloor\]
\[ R_{{{6}}}^{[5]}= \pi_{{{[1,2,3,4,5]}}}^{[1,2,3,4,5]}\lfloor
\pi_{{{[2,5]}}}^{[1,2]}   R_{{{1}}}^{[2]}, \pi_{{{[1,3,4]}}}^{[1,2,3]}
   R_{{{1}}}^{[3]}\rfloor\]
\[ R_{{{7}}}^{[5]}= \pi_{{{[1,2,3,4,5]}}}^{[1,2,3,4,5]}\lfloor
\pi_{{{[1,5]}}}^{[1,2]}   R_{{{1}}}^{[2]}, \pi_{{{[2,3,4]}}}^{[1,2,3]}
   R_{{{1}}}^{[3]}\rfloor\]
\[ R_{{{8}}}^{[5]}= \pi_{{{[1,2,3,4,5]}}}^{[1,2,3,4,5]}\lfloor
\pi_{{{[2,4]}}}^{[1,2]}   R_{{{1}}}^{[2]}, \pi_{{{[1,3,5]}}}^{[1,2,3]}
   R_{{{1}}}^{[3]}\rfloor\]

(e) The top canonical generating system for $\mathcal R_6$ is
given by
$R_{{{1}}}^{[2]}$
and $R_{{{1}}}^{[6]},\ldots, R_{{{15}}}^{[6]}$; where:
\[ R_{{{1}}}^{[6]}= \pi_{{{[1,1,2,3,4,4]}}}^{[1,2,3,4,5,6]}\lfloor
 \pi_{{{[3]}}}^{[1]}   R_{{{1}}}^{[1]}, \pi_{{{[1,2,4,5,6]}}}^{[1,2,3
,4,5]}   R_{{{5}}}^{[5]}\rfloor\]
\[ R_{{{2}}}^{[6]}= \pi_{{{[1,1,2,3,4,4]}}}^{[1,2,3,4,5,6]}\lfloor
 \pi_{{{[4]}}}^{[1]}   R_{{{1}}}^{[1]}, \pi_{{{[1,2,3,5,6]}}}^{[1,2,3
,4,5]}   R_{{{5}}}^{[5]}\rfloor\]
\[ R_{{{3}}}^{[6]}= \pi_{{{[1,2,3,4,4]}}}^{[2,3,4,5,6]}\lfloor
\pi_{{{[2]}}}^{[1]}   R_{{{1}}}^{[1]}, \pi_{{{[1,3,4,5,6]}}}^{[1,2,3,4,5
]}   R_{{{5}}}^{[5]}\rfloor\]
\[ R_{{{4}}}^{[6]}= \pi_{{{[1,1,2,3,4]}}}^{[1,2,3,4,5]}\lfloor
\pi_{{{[5]}}}^{[1]}   R_{{{1}}}^{[1]}, \pi_{{{[1,2,3,4,6]}}}^{[1,2,3,4,5
]}   R_{{{6}}}^{[5]}\rfloor\]
\[ R_{{{5}}}^{[6]}= \pi_{{{[1,2,3,4,5]}}}^{[1,2,3,4,5]}\lfloor
\pi_{{{[4]}}}^{[1]}   R_{{{1}}}^{[1]}, \pi_{{{[1,2,3,5]}}}^{[1,2,3,4]}
   R_{{{3}}}^{[5]}\rfloor\]
\[ R_{{{6}}}^{[6]}= \pi_{{{[1,2,3,4,5]}}}^{[1,2,3,4,5]}\lfloor
\pi_{{{[2]}}}^{[1]}   R_{{{1}}}^{[1]}, \pi_{{{[1,3,4,5]}}}^{[1,2,3,4]}
   R_{{{4}}}^{[5]}\rfloor\]
\[ R_{{{7}}}^{[6]}= \pi_{{{[1,2,3,4,5]}}}^{[1,2,3,4,5]}\lfloor
\pi_{{{[1,2,6]}}}^{[1,2,3]}   R_{{{1}}}^{[3]}, \pi_{{{[3,4,5]}}}^{[1,2,3
]}   R_{{{1}}}^{[3]}\rfloor\]
\[ R_{{{8}}}^{[6]}= \pi_{{{[1,2,3,4,5]}}}^{[1,3,4,5,6]}\lfloor
\pi_{{{[1,3,4]}}}^{[1,2,3]}   R_{{{1}}}^{[3]}, \pi_{{{[2,5,6]}}}^{[1,2,3
]}   R_{{{1}}}^{[3]}\rfloor\]
\[ R_{{{9}}}^{[6]}= \pi_{{{[1,2,3,4,5,6]}}}^{[1,2,3,4,5,6]}\lfloor
 \pi_{{{[1,4,5]}}}^{[1,2,3]}   R_{{{1}}}^{[3]}, \pi_{{{[2,3,6]}}}^{[1
,2,3]}   R_{{{1}}}^{[3]}\rfloor\]
\[ R_{{{10}}}^{[6]}= \pi_{{{[1,2,3,4,5,6]}}}^{[1,2,3,4,5,6]}
\lfloor \pi_{{{[1,4,6]}}}^{[1,2,3]}   R_{{{1}}}^{[3]}, \pi_{{{[2,3,5]}}}^{[
1,2,3]}   R_{{{1}}}^{[3]}\rfloor\]
\[ R_{{{11}}}^{[6]}= \pi_{{{[1,2,3,4,5,6]}}}^{[1,2,3,4,5,6]}
\lfloor \pi_{{{[1,3,6]}}}^{[1,2,3]}   R_{{{1}}}^{[3]}, \pi_{{{[2,4,5]}}}^{[
1,2,3]}   R_{{{1}}}^{[3]}\rfloor\]
\[ R_{{{12}}}^{[6]}= \pi_{{{[1,2,3,4,5,6]}}}^{[1,2,3,4,5,6]}
\lfloor \pi_{{{[1,3,5]}}}^{[1,2,3]}   R_{{{1}}}^{[3]}, \pi_{{{[2,4,6]}}}^{[
1,2,3]}   R_{{{1}}}^{[3]}\rfloor\]
\[ R_{{{13}}}^{[6]}= \pi_{{{[1,2]}}}^{[1,2]}\lfloor   \pi_{{{[2]}}}^{
[1]}   R_{{{1}}}^{[1]}, \pi_{{{[1,3]}}}^{[1,2]}   R_{{{2}}}^{[5]}
\rfloor\]
\[ R_{{{14}}}^{[6]}= \pi_{{{[1,2,3,4,4]}}}^{[1,2,3,4,5]}\lfloor  {
}\pi_{{{[3]}}}^{[1]}   R_{{{1}}}^{[1]}, \pi_{{{[1,2,4,5]}}}^{[1,2,3,4]
}   R_{{{3}}}^{[5]}\rfloor\]
\[ R_{{{15}}}^{[6]}= \pi_{{{[1,1,2,3,4]}}}^{[1,2,3,4,5]}\lfloor  {
}\pi_{{{[3]}}}^{[1]}   R_{{{1}}}^{[1]}, \pi_{{{[1,2,4,5]}}}^{[1,2,3,4]
}   R_{{{4}}}^{[5]}\rfloor\]
Remark:
\begin{equation}
R_{{1}}^{[2]}=\pi _{{[1,1,2]}}^{[1,2,3]}\lfloor \pi _{{[3]}}^{[1]
} R_{{1}}^{[1]},\pi _{{[1,2]}}^{[1,2]} R_{{1}}^{[5]}\rfloor.
\plabel{eq:gen26}
\end{equation}

(f) The top canonical generating system for $\mathcal R_7$ is
given by
$R_{{{1}}}^{[3]}$
and $R_{{{1}}}^{[7]},\ldots, R_{{{28}}}^{[7]}$; where:
\[ R_{{{1}}}^{[7]}= \pi_{{{[1,2,3,4]}}}^{[2,3,4,5]}\lfloor   \pi_{{{[
4]}}}^{[1]}   R_{{{1}}}^{[1]}, \pi_{{{[1,2,3,5]}}}^{[1,2,3,4]}
 R_{{{1}}}^{[6]}\rfloor\]
\[ R_{{{2}}}^{[7]}= \pi_{{{[1,2,3,4]}}}^{[1,2,3,4]}\lfloor   \pi_{{{[
2]}}}^{[1]}   R_{{{1}}}^{[1]}, \pi_{{{[1,3,4,5]}}}^{[1,2,3,4]}
 R_{{{2}}}^{[6]}\rfloor\]
\[ R_{{{3}}}^{[7]}= \pi_{{{[1,1,2,3,3,4]}}}^{[1,2,3,4,5,6]}\lfloor
 \pi_{{{[3]}}}^{[1]}   R_{{{1}}}^{[1]}, \pi_{{{[1,2,4,5,6,7]}}}^{[1,2
,3,4,5,6]}   R_{{{9}}}^{[6]}\rfloor\]
\[ R_{{{4}}}^{[7]}= \pi_{{{[1,2,2,3,4,4]}}}^{[1,3,4,5,6,7]}\lfloor
 \pi_{{{[5]}}}^{[1]}   R_{{{1}}}^{[1]}, \pi_{{{[1,2,3,4,6,7]}}}^{[1,2
,3,4,5,6]}   R_{{{11}}}^{[6]}\rfloor\]
\[ R_{{{5}}}^{[7]}= \pi_{{{[1,2,3,3,4,5]}}}^{[1,3,4,5,6,7]}\lfloor
 \pi_{{{[6]}}}^{[1]}   R_{{{1}}}^{[1]}, \pi_{{{[1,2,3,4,5,7]}}}^{[1,2
,3,4,5,6]}   R_{{{12}}}^{[6]}\rfloor\]
\[ R_{{{6}}}^{[7]}= \pi_{{{[1,2,3,3,4,5]}}}^{[1,2,3,4,5,6]}\lfloor
 \pi_{{{[2]}}}^{[1]}   R_{{{1}}}^{[1]}, \pi_{{{[1,3,4,5,6,7]}}}^{[1,2
,3,4,5,6]}   R_{{{11}}}^{[6]}\rfloor\]
\[ R_{{{7}}}^{[7]}= \pi_{{{[1,2,3,4,5]}}}^{[1,2,3,4,5]}\lfloor
\pi_{{{[1,2,6]}}}^{[1,2,3]}   R_{{{1}}}^{[3]}, \pi_{{{[3,4,5]}}}^{[1,2,3
]}   R_{{{1}}}^{[4]}\rfloor\]
\[ R_{{{8}}}^{[7]}= \pi_{{{[1,2,3,4,5]}}}^{[1,3,4,5,6]}\lfloor
\pi_{{{[2,5,6]}}}^{[1,2,3]}   R_{{{1}}}^{[3]}, \pi_{{{[1,3,4]}}}^{[1,2,3
]}   R_{{{2}}}^{[4]}\rfloor\]
\[ R_{{{9}}}^{[7]}= \pi_{{{[1,2,3,4,5]}}}^{[2,3,4,5,6]}\lfloor
\pi_{{{[1,5,6]}}}^{[1,2,3]}   R_{{{1}}}^{[3]}, \pi_{{{[2,3,4]}}}^{[1,2,3
]}   R_{{{1}}}^{[4]}\rfloor\]
\[ R_{{{10}}}^{[7]}= \pi_{{{[1,2,3,4,5]}}}^{[1,2,3,4,5]}\lfloor  {
}\pi_{{{[1,2,6]}}}^{[1,2,3]}   R_{{{1}}}^{[3]}, \pi_{{{[3,4,5]}}}^{[1,
2,3]}   R_{{{2}}}^{[4]}\rfloor\]
\[ R_{{{11}}}^{[7]}= \pi_{{{[1,2,3,4,5,6]}}}^{[1,2,3,4,5,6]}
\lfloor \pi_{{{[2,4,5]}}}^{[1,2,3]}   R_{{{1}}}^{[3]}, \pi_{{{[1,3,6]}}}^{[
1,2,3]}   R_{{{1}}}^{[4]}\rfloor\]
\[ R_{{{12}}}^{[7]}= \pi_{{{[1,2,3,4,5,6]}}}^{[1,2,3,4,5,6]}
\lfloor \pi_{{{[2,4,6]}}}^{[1,2,3]}   R_{{{1}}}^{[3]}, \pi_{{{[1,3,5]}}}^{[
1,2,3]}   R_{{{1}}}^{[4]}\rfloor\]
\[ R_{{{13}}}^{[7]}= \pi_{{{[1,2,3,4,5,6]}}}^{[1,2,3,4,5,6]}
\lfloor \pi_{{{[1,4,5]}}}^{[1,2,3]}   R_{{{1}}}^{[3]}, \pi_{{{[2,3,6]}}}^{[
1,2,3]}   R_{{{1}}}^{[4]}\rfloor\]
\[ R_{{{14}}}^{[7]}= \pi_{{{[1,2,3,4,5,6]}}}^{[1,2,3,4,5,6]}
\lfloor \pi_{{{[1,4,6]}}}^{[1,2,3]}   R_{{{1}}}^{[3]}, \pi_{{{[2,3,5]}}}^{[
1,2,3]}   R_{{{1}}}^{[4]}\rfloor\]
\[ R_{{{15}}}^{[7]}= \pi_{{{[1,2,3,4,5,6]}}}^{[1,2,3,4,5,6]}
\lfloor \pi_{{{[2,3,6]}}}^{[1,2,3]}   R_{{{1}}}^{[3]}, \pi_{{{[1,4,5]}}}^{[
1,2,3]}   R_{{{2}}}^{[4]}\rfloor\]
\[ R_{{{16}}}^{[7]}= \pi_{{{[1,2,3,4,5,6]}}}^{[1,2,3,4,5,6]}
\lfloor \pi_{{{[2,3,5]}}}^{[1,2,3]}   R_{{{1}}}^{[3]}, \pi_{{{[1,4,6]}}}^{[
1,2,3]}   R_{{{2}}}^{[4]}\rfloor\]
\[ R_{{{17}}}^{[7]}= \pi_{{{[1,2,3,4,5,6]}}}^{[1,2,3,4,5,6]}
\lfloor \pi_{{{[1,3,6]}}}^{[1,2,3]}   R_{{{1}}}^{[3]}, \pi_{{{[2,4,5]}}}^{[
1,2,3]}   R_{{{2}}}^{[4]}\rfloor\]
\[ R_{{{18}}}^{[7]}= \pi_{{{[1,2,3,4,5,6]}}}^{[1,2,3,4,5,6]}
\lfloor \pi_{{{[1,3,5]}}}^{[1,2,3]}   R_{{{1}}}^{[3]}, \pi_{{{[2,4,6]}}}^{[
1,2,3]}   R_{{{2}}}^{[4]}\rfloor\]
\[ R_{{{19}}}^{[7]}= \pi_{{{[1,2,2,2]}}}^{[1,2,3,4]}\lfloor   \pi_{{{
[1]}}}^{[1]}   R_{{{1}}}^{[1]}, \pi_{{{[2,3,4,5]}}}^{[1,2,3,4]}
   R_{{{1}}}^{[6]}\rfloor\]
\[ R_{{{20}}}^{[7]}= \pi_{{{[1,1,1,2]}}}^{[2,3,4,5]}\lfloor   \pi_{{{
[5]}}}^{[1]}   R_{{{1}}}^{[1]}, \pi_{{{[1,2,3,4]}}}^{[1,2,3,4]}
   R_{{{1}}}^{[6]}\rfloor\]
\[ R_{{{21}}}^{[7]}= \pi_{{{[1,2,2,3,4,5,5]}}}^{[1,2,3,4,5,6,7]}
\lfloor   \pi_{{{[1,5]}}}^{[1,2]}   R_{{{1}}}^{[2]}, \pi_{{{[2,3,4,6,7
]}}}^{[1,2,3,4,5]}   R_{{{5}}}^{[5]}\rfloor\]
\[ R_{{{22}}}^{[7]}= \pi_{{{[1,2,2,3,4,4,5]}}}^{[1,2,3,4,5,6,7]}
\lfloor   \pi_{{{[1,7]}}}^{[1,2]}   R_{{{1}}}^{[2]}, \pi_{{{[2,3,4,5,6
]}}}^{[1,2,3,4,5]}   R_{{{5}}}^{[5]}\rfloor\]
\[ R_{{{23}}}^{[7]}= \pi_{{{[1,1,2,3,4,5,5]}}}^{[1,2,3,4,5,6,7]}
\lfloor   \pi_{{{[3,5]}}}^{[1,2]}   R_{{{1}}}^{[2]}, \pi_{{{[1,2,4,6,7
]}}}^{[1,2,3,4,5]}   R_{{{5}}}^{[5]}\rfloor\]
\[ R_{{{24}}}^{[7]}= \pi_{{{[1,1,2,3,4,4,5]}}}^{[1,2,3,4,5,6,7]}
\lfloor   \pi_{{{[3,7]}}}^{[1,2]}   R_{{{1}}}^{[2]}, \pi_{{{[1,2,4,5,6
]}}}^{[1,2,3,4,5]}   R_{{{5}}}^{[5]}\rfloor\]
\[ R_{{{25}}}^{[7]}= \pi_{{{[1,2,3,4,5,6]}}}^{[1,2,3,4,5,6]}
\lfloor \pi_{{{[2,3,5]}}}^{[1,2,3]}   R_{{{1}}}^{[3]}, \pi_{{{[1,4,6]}}}^{[
1,2,3]}   R_{{{1}}}^{[4]}\rfloor\]
\[ R_{{{26}}}^{[7]}= \pi_{{{[1,2,3,4,5,6]}}}^{[1,2,3,4,5,6]}
\lfloor \pi_{{{[1,3,4]}}}^{[1,2,3]}   R_{{{1}}}^{[3]}, \pi_{{{[2,5,6]}}}^{[
1,2,3]}   R_{{{1}}}^{[4]}\rfloor\]
\[ R_{{{27}}}^{[7]}= \pi_{{{[1,2,3,4,5,6]}}}^{[1,2,3,4,5,6]}
\lfloor \pi_{{{[2,4,5]}}}^{[1,2,3]}   R_{{{1}}}^{[3]}, \pi_{{{[1,3,6]}}}^{[
1,2,3]}   R_{{{2}}}^{[4]}\rfloor\]
\[ R_{{{28}}}^{[7]}= \pi_{{{[1,2,3,4,5,6]}}}^{[1,2,3,4,5,6]}
\lfloor \pi_{{{[3,4,6]}}}^{[1,2,3]}   R_{{{1}}}^{[3]}, \pi_{{{[1,2,5]}}}^{[
1,2,3]}   R_{{{2}}}^{[4]}\rfloor\]

(g) The top canonical generating system for $\mathcal R_8$ is
given by
$R_{{{1}}}^{[4]},R_{{{2}}}^{[4]} $
and $R_{{{1}}}^{[8]},\ldots, R_{{{31}}}^{[8]}$; where:
\[ R_{{{1}}}^{[8]}= \pi_{{{[1,1,2,2]}}}^{[2,3,4,5]}\lfloor   \pi_{{{[
4,5,6]}}}^{[1,2,3]}   R_{{{1}}}^{[4]}, \pi_{{{[1,2,3]}}}^{[1,2,3]}
   R_{{{2}}}^{[4]}\rfloor\]
\[ R_{{{2}}}^{[8]}= \pi_{{{[1,2,3,4]}}}^{[1,2,3,4]}\lfloor   \pi_{{{[
4]}}}^{[1]}   R_{{{1}}}^{[1]}, \pi_{{{[1,2,3,5]}}}^{[1,2,3,4]}
 R_{{{3}}}^{[7]}\rfloor\]
\[ R_{{{3}}}^{[8]}= \pi_{{{[1,2,3,4]}}}^{[2,3,4,5]}\lfloor   \pi_{{{[
2]}}}^{[1]}   R_{{{1}}}^{[1]}, \pi_{{{[1,3,4,5]}}}^{[1,2,3,4]}
 R_{{{4}}}^{[7]}\rfloor\]
\[ R_{{{4}}}^{[8]}= \pi_{{{[1,2,2,3,4]}}}^{[2,3,4,5,6]}\lfloor
\pi_{{{[2]}}}^{[1]}   R_{{{1}}}^{[1]}, \pi_{{{[1,3,4,5,6]}}}^{[1,2,3,4,5
]}   R_{{{5}}}^{[7]}\rfloor\]
\[ R_{{{5}}}^{[8]}= \pi_{{{[1,2,3,3,4]}}}^{[1,2,3,4,5]}\lfloor
\pi_{{{[5]}}}^{[1]}   R_{{{1}}}^{[1]}, \pi_{{{[1,2,3,4,6]}}}^{[1,2,3,4,5
]}   R_{{{6}}}^{[7]}\rfloor\]
\[ R_{{{6}}}^{[8]}= \pi_{{{[1,1,2,3,4,4]}}}^{[1,2,3,4,5,6]}\lfloor
 \pi_{{{[3]}}}^{[1]}   R_{{{1}}}^{[1]}, \pi_{{{[1,2,4,5,6]}}}^{[1,2,3
,4,5]}   R_{{{21}}}^{[7]}\rfloor\]
\[ R_{{{7}}}^{[8]}= \pi_{{{[1,1,2,3,4,4]}}}^{[1,2,3,4,5,6]}\lfloor
 \pi_{{{[4]}}}^{[1]}   R_{{{1}}}^{[1]}, \pi_{{{[1,2,3,5,6]}}}^{[1,2,3
,4,5]}   R_{{{21}}}^{[7]}\rfloor\]
\[ R_{{{8}}}^{[8]}= \pi_{{{[1,2,2,3,4,5,6,7]}}}^{[1,2,3,4,5,6,7,8]}
\lfloor   \pi_{{{[1,5,8]}}}^{[1,2,3]}   R_{{{1}}}^{[3]}, \pi_{{{[2,3,4
,6,7]}}}^{[1,2,3,4,5]}   R_{{{6}}}^{[5]}\rfloor\]
\[ R_{{{9}}}^{[8]}= \pi_{{{[1,2,2,3,4,5,6,7]}}}^{[1,2,3,4,5,6,7,8]}
\lfloor   \pi_{{{[1,5,8]}}}^{[1,2,3]}   R_{{{1}}}^{[3]}, \pi_{{{[2,3,4
,6,7]}}}^{[1,2,3,4,5]}   R_{{{5}}}^{[5]}\rfloor\]
\[ R_{{{10}}}^{[8]}= \pi_{{{[1,2,2,3,4,5,6,7]}}}^{[1,2,3,4,5,6,7,8]
}\lfloor   \pi_{{{[1,5,7]}}}^{[1,2,3]}   R_{{{1}}}^{[3]}, \pi_{{{[2,3,
4,6,8]}}}^{[1,2,3,4,5]}   R_{{{6}}}^{[5]}\rfloor\]
\[ R_{{{11}}}^{[8]}= \pi_{{{[1,2,2,3,4,5,6,7]}}}^{[1,2,3,4,5,6,7,8]
}\lfloor   \pi_{{{[1,5,7]}}}^{[1,2,3]}   R_{{{1}}}^{[3]}, \pi_{{{[2,3,
4,6,8]}}}^{[1,2,3,4,5]}   R_{{{5}}}^{[5]}\rfloor\]
\[ R_{{{12}}}^{[8]}= \pi_{{{[1,2,2,3,4,5,6,7]}}}^{[1,2,3,4,5,6,7,8]
}\lfloor   \pi_{{{[1,5,6]}}}^{[1,2,3]}   R_{{{1}}}^{[3]}, \pi_{{{[2,3,
4,7,8]}}}^{[1,2,3,4,5]}   R_{{{6}}}^{[5]}\rfloor\]
\[ R_{{{13}}}^{[8]}= \pi_{{{[1,2,2,3,4,5,6,7]}}}^{[1,2,3,4,5,6,7,8]
}\lfloor   \pi_{{{[1,5,6]}}}^{[1,2,3]}   R_{{{1}}}^{[3]}, \pi_{{{[2,3,
4,7,8]}}}^{[1,2,3,4,5]}   R_{{{5}}}^{[5]}\rfloor\]
\[ R_{{{14}}}^{[8]}= \pi_{{{[1,2,3,4,5,6,7,7]}}}^{[1,2,3,4,5,6,7,8]
}\lfloor   \pi_{{{[3,4,6]}}}^{[1,2,3]}   R_{{{1}}}^{[3]}, \pi_{{{[1,2,
5,7,8]}}}^{[1,2,3,4,5]}   R_{{{6}}}^{[5]}\rfloor\]
\[ R_{{{15}}}^{[8]}= \pi_{{{[1,2,3,4,5,6,7,7]}}}^{[1,2,3,4,5,6,7,8]
}\lfloor   \pi_{{{[2,4,6]}}}^{[1,2,3]}   R_{{{1}}}^{[3]}, \pi_{{{[1,3,
5,7,8]}}}^{[1,2,3,4,5]}   R_{{{6}}}^{[5]}\rfloor\]
\[ R_{{{16}}}^{[8]}= \pi_{{{[1,2,3,4,5,6,6,7]}}}^{[1,2,3,4,5,6,7,8]
}\lfloor   \pi_{{{[3,4,8]}}}^{[1,2,3]}   R_{{{1}}}^{[3]}, \pi_{{{[1,2,
5,6,7]}}}^{[1,2,3,4,5]}   R_{{{6}}}^{[5]}\rfloor\]
\[ R_{{{17}}}^{[8]}= \pi_{{{[1,2,3,4,5,6,6,7]}}}^{[1,2,3,4,5,6,7,8]
}\lfloor   \pi_{{{[2,4,8]}}}^{[1,2,3]}   R_{{{1}}}^{[3]}, \pi_{{{[1,3,
5,6,7]}}}^{[1,2,3,4,5]}   R_{{{6}}}^{[5]}\rfloor\]
\[ R_{{{18}}}^{[8]}= \pi_{{{[1,2,3,4,5,6,7,7]}}}^{[1,2,3,4,5,6,7,8]
}\lfloor   \pi_{{{[3,4,6]}}}^{[1,2,3]}   R_{{{1}}}^{[3]}, \pi_{{{[1,2,
5,7,8]}}}^{[1,2,3,4,5]}   R_{{{5}}}^{[5]}\rfloor\]
\[ R_{{{19}}}^{[8]}= \pi_{{{[1,2,3,4,5,6,7,7]}}}^{[1,2,3,4,5,6,7,8]
}\lfloor   \pi_{{{[2,4,6]}}}^{[1,2,3]}   R_{{{1}}}^{[3]}, \pi_{{{[1,3,
5,7,8]}}}^{[1,2,3,4,5]}   R_{{{5}}}^{[5]}\rfloor\]
\[ R_{{{20}}}^{[8]}= \pi_{{{[1,2,3,4,5,6,6,7]}}}^{[1,2,3,4,5,6,7,8]
}\lfloor   \pi_{{{[3,4,8]}}}^{[1,2,3]}   R_{{{1}}}^{[3]}, \pi_{{{[1,2,
5,6,7]}}}^{[1,2,3,4,5]}   R_{{{5}}}^{[5]}\rfloor\]
\[ R_{{{21}}}^{[8]}= \pi_{{{[1,2,3,4,5,6,6,7]}}}^{[1,2,3,4,5,6,7,8]
}\lfloor   \pi_{{{[2,4,8]}}}^{[1,2,3]}   R_{{{1}}}^{[3]}, \pi_{{{[1,3,
5,6,7]}}}^{[1,2,3,4,5]}   R_{{{5}}}^{[5]}\rfloor\]
\[ R_{{{22}}}^{[8]}= \pi_{{{[1,1,2,3,4,5,6,7]}}}^{[1,2,3,4,5,6,7,8]
}\lfloor   \pi_{{{[3,5,6]}}}^{[1,2,3]}   R_{{{1}}}^{[3]}, \pi_{{{[1,2,
4,7,8]}}}^{[1,2,3,4,5]}   R_{{{6}}}^{[5]}\rfloor\]
\[ R_{{{23}}}^{[8]}= \pi_{{{[1,1,2,3,4,5,6,7]}}}^{[1,2,3,4,5,6,7,8]
}\lfloor   \pi_{{{[3,5,6]}}}^{[1,2,3]}   R_{{{1}}}^{[3]}, \pi_{{{[1,2,
4,7,8]}}}^{[1,2,3,4,5]}   R_{{{5}}}^{[5]}\rfloor\]
\[ R_{{{24}}}^{[8]}= \pi_{{{[1,1,2,3,4,5,6,7]}}}^{[1,2,3,4,5,6,7,8]
}\lfloor   \pi_{{{[3,5,7]}}}^{[1,2,3]}   R_{{{1}}}^{[3]}, \pi_{{{[1,2,
4,6,8]}}}^{[1,2,3,4,5]}   R_{{{6}}}^{[5]}\rfloor\]
\[ R_{{{25}}}^{[8]}= \pi_{{{[1,1,2,3,4,5,6,7]}}}^{[1,2,3,4,5,6,7,8]
}\lfloor   \pi_{{{[3,5,7]}}}^{[1,2,3]}   R_{{{1}}}^{[3]}, \pi_{{{[1,2,
4,6,8]}}}^{[1,2,3,4,5]}   R_{{{5}}}^{[5]}\rfloor\]
\[ R_{{{26}}}^{[8]}= \pi_{{{[1,1,2,3,4,5,6,7]}}}^{[1,2,3,4,5,6,7,8]
}\lfloor   \pi_{{{[3,5,8]}}}^{[1,2,3]}   R_{{{1}}}^{[3]}, \pi_{{{[1,2,
4,6,7]}}}^{[1,2,3,4,5]}   R_{{{6}}}^{[5]}\rfloor\]
\[ R_{{{27}}}^{[8]}= \pi_{{{[1,1,2,3,4,5,6,7]}}}^{[1,2,3,4,5,6,7,8]
}\lfloor   \pi_{{{[3,5,8]}}}^{[1,2,3]}   R_{{{1}}}^{[3]}, \pi_{{{[1,2,
4,6,7]}}}^{[1,2,3,4,5]}   R_{{{5}}}^{[5]}\rfloor\]
\[ R_{{{28}}}^{[8]}= \pi_{{{[1,2,3,4,5,6,6,7]}}}^{[1,2,3,4,5,6,7,8]
}\lfloor   \pi_{{{[1,4,8]}}}^{[1,2,3]}   R_{{{1}}}^{[3]}, \pi_{{{[2,3,
5,6,7]}}}^{[1,2,3,4,5]}   R_{{{6}}}^{[5]}\rfloor\]
\[ R_{{{29}}}^{[8]}= \pi_{{{[1,2,3,4,5,6,6,7]}}}^{[1,2,3,4,5,6,7,8]
}\lfloor   \pi_{{{[1,4,8]}}}^{[1,2,3]}   R_{{{1}}}^{[3]}, \pi_{{{[2,3,
5,6,7]}}}^{[1,2,3,4,5]}   R_{{{5}}}^{[5]}\rfloor\]
\[ R_{{{30}}}^{[8]}= \pi_{{{[1,2,3,4,5,6,7,7]}}}^{[1,2,3,4,5,6,7,8]
}\lfloor   \pi_{{{[1,4,6]}}}^{[1,2,3]}   R_{{{1}}}^{[3]}, \pi_{{{[2,3,
5,7,8]}}}^{[1,2,3,4,5]}   R_{{{6}}}^{[5]}\rfloor\]
\[ R_{{{31}}}^{[8]}= \pi_{{{[1,2,3,4,5,6,7,7]}}}^{[1,2,3,4,5,6,7,8]
}\lfloor   \pi_{{{[1,4,6]}}}^{[1,2,3]}   R_{{{1}}}^{[3]}, \pi_{{{[2,3,
5,7,8]}}}^{[1,2,3,4,5]}   R_{{{5}}}^{[5]}\rfloor\]

(h) The top canonical generating system for $\mathcal R_9$ is
given by
$R_{{{1}}}^{[5]},\ldots,R_{{{8}}}^{[5]} $
and $R_{{{1}}}^{[9]},\ldots, R_{{{81}}}^{[9]}$; where:
\[ R_{{{1}}}^{[9]}= \pi_{{{[1,1,1,2,3]}}}^{[1,2,3,4,5]}\lfloor
\pi_{{{[4]}}}^{[1]}   R_{{{1}}}^{[1]}, \pi_{{{[1,2,3,5]}}}^{[1,2,3,4]}
   R_{{{2}}}^{[8]}\rfloor\]
\[ R_{{{2}}}^{[9]}= \pi_{{{[1,2,3,3,4,4,5,5]}}}^{[1,2,3,4,5,6,7,8]}
\lfloor   \pi_{{{[1,2,9]}}}^{[1,2,3]}   R_{{{1}}}^{[3]}, \pi_{{{[3,4,5
,6,7,8]}}}^{[1,2,3,4,5,6]}   R_{{{9}}}^{[6]}\rfloor\]
\[ R_{{{3}}}^{[9]}= \pi_{{{[1,1,2,2,3,3,4,5]}}}^{[1,2,4,5,6,7,8,9]}
\lfloor   \pi_{{{[3,8,9]}}}^{[1,2,3]}   R_{{{1}}}^{[3]}, \pi_{{{[1,2,4
,5,6,7]}}}^{[1,2,3,4,5,6]}   R_{{{9}}}^{[6]}\rfloor\]
\[ R_{{{4}}}^{[9]}= \pi_{{{[1,2,3,4,4,5,6,7,7]}}}^{[1,2,3,4,5,6,7,8
,9]}\lfloor   \pi_{{{[1,6,7]}}}^{[1,2,3]}   R_{{{1}}}^{[3]}, \pi_{{{[2
,3,4,5,8,9]}}}^{[1,2,3,4,5,6]}   R_{{{11}}}^{[6]}\rfloor\]
\[ R_{{{5}}}^{[9]}= \pi_{{{[1,2,3,4,4,5,6,6,7]}}}^{[1,2,3,4,5,6,7,8
,9]}\lfloor   \pi_{{{[1,6,9]}}}^{[1,2,3]}   R_{{{1}}}^{[3]}, \pi_{{{[2
,3,4,5,7,8]}}}^{[1,2,3,4,5,6]}   R_{{{11}}}^{[6]}\rfloor\]
\[ R_{{{6}}}^{[9]}= \pi_{{{[1,2,3,4,4,5,6,6,7]}}}^{[1,2,3,4,5,6,7,8
,9]}\lfloor   \pi_{{{[2,6,9]}}}^{[1,2,3]}   R_{{{1}}}^{[3]}, \pi_{{{[1
,3,4,5,7,8]}}}^{[1,2,3,4,5,6]}   R_{{{11}}}^{[6]}\rfloor\]
\[ R_{{{7}}}^{[9]}= \pi_{{{[1,2,3,4,4,5,6,7,7]}}}^{[1,2,3,4,5,6,7,8
,9]}\lfloor   \pi_{{{[2,6,7]}}}^{[1,2,3]}   R_{{{1}}}^{[3]}, \pi_{{{[1
,3,4,5,8,9]}}}^{[1,2,3,4,5,6]}   R_{{{11}}}^{[6]}\rfloor\]
\[ R_{{{8}}}^{[9]}= \pi_{{{[1,2,3,4,4,5,6,6,7]}}}^{[1,2,3,4,5,6,7,8
,9]}\lfloor   \pi_{{{[1,6,9]}}}^{[1,2,3]}   R_{{{1}}}^{[3]}, \pi_{{{[2
,3,4,5,7,8]}}}^{[1,2,3,4,5,6]}   R_{{{9}}}^{[6]}\rfloor\]
\[ R_{{{9}}}^{[9]}= \pi_{{{[1,2,3,4,4,5,6,7,7]}}}^{[1,2,3,4,5,6,7,8
,9]}\lfloor   \pi_{{{[1,6,7]}}}^{[1,2,3]}   R_{{{1}}}^{[3]}, \pi_{{{[2
,3,4,5,8,9]}}}^{[1,2,3,4,5,6]}   R_{{{9}}}^{[6]}\rfloor\]
\[ R_{{{10}}}^{[9]}= \pi_{{{[1,2,3,4,4,5,6,7,7]}}}^{[1,2,3,4,5,6,7,
8,9]}\lfloor   \pi_{{{[2,6,7]}}}^{[1,2,3]}   R_{{{1}}}^{[3]}, \pi_{{{[
1,3,4,5,8,9]}}}^{[1,2,3,4,5,6]}   R_{{{9}}}^{[6]}\rfloor\]
\[ R_{{{11}}}^{[9]}= \pi_{{{[1,2,3,4,4,5,6,6,7]}}}^{[1,2,3,4,5,6,7,
8,9]}\lfloor   \pi_{{{[2,6,9]}}}^{[1,2,3]}   R_{{{1}}}^{[3]}, \pi_{{{[
1,3,4,5,7,8]}}}^{[1,2,3,4,5,6]}   R_{{{9}}}^{[6]}\rfloor\]
\[ R_{{{12}}}^{[9]}= \pi_{{{[1,2,3,4,4,5,6,7,7]}}}^{[1,2,3,4,5,6,7,
8,9]}\lfloor   \pi_{{{[3,6,7]}}}^{[1,2,3]}   R_{{{1}}}^{[3]}, \pi_{{{[
1,2,4,5,8,9]}}}^{[1,2,3,4,5,6]}   R_{{{11}}}^{[6]}\rfloor\]
\[ R_{{{13}}}^{[9]}= \pi_{{{[1,2,3,4,4,5,6,6,7]}}}^{[1,2,3,4,5,6,7,
8,9]}\lfloor   \pi_{{{[3,6,9]}}}^{[1,2,3]}   R_{{{1}}}^{[3]}, \pi_{{{[
1,2,4,5,7,8]}}}^{[1,2,3,4,5,6]}   R_{{{11}}}^{[6]}\rfloor\]
\[ R_{{{14}}}^{[9]}= \pi_{{{[1,2,3,4,4,5,6,6,7]}}}^{[1,2,3,4,5,6,7,
8,9]}\lfloor   \pi_{{{[3,6,9]}}}^{[1,2,3]}   R_{{{1}}}^{[3]}, \pi_{{{[
1,2,4,5,7,8]}}}^{[1,2,3,4,5,6]}   R_{{{9}}}^{[6]}\rfloor\]
\[ R_{{{15}}}^{[9]}= \pi_{{{[1,2,3,4,4,5,6,7,7]}}}^{[1,2,3,4,5,6,7,
8,9]}\lfloor   \pi_{{{[3,6,7]}}}^{[1,2,3]}   R_{{{1}}}^{[3]}, \pi_{{{[
1,2,4,5,8,9]}}}^{[1,2,3,4,5,6]}   R_{{{9}}}^{[6]}\rfloor\]
\[ R_{{{16}}}^{[9]}= \pi_{{{[1,1,2,3,4,4,5,6,7]}}}^{[1,2,3,4,5,6,7,
8,9]}\lfloor   \pi_{{{[3,4,8]}}}^{[1,2,3]}   R_{{{1}}}^{[3]}, \pi_{{{[
1,2,5,6,7,9]}}}^{[1,2,3,4,5,6]}   R_{{{9}}}^{[6]}\rfloor\]
\[ R_{{{17}}}^{[9]}= \pi_{{{[1,1,2,3,4,4,5,6,7]}}}^{[1,2,3,4,5,6,7,
8,9]}\lfloor   \pi_{{{[3,4,9]}}}^{[1,2,3]}   R_{{{1}}}^{[3]}, \pi_{{{[
1,2,5,6,7,8]}}}^{[1,2,3,4,5,6]}   R_{{{10}}}^{[6]}\rfloor\]
\[ R_{{{18}}}^{[9]}= \pi_{{{[1,1,2,3,4,4,5,6,7]}}}^{[1,2,3,4,5,6,7,
8,9]}\lfloor   \pi_{{{[3,4,8]}}}^{[1,2,3]}   R_{{{1}}}^{[3]}, \pi_{{{[
1,2,5,6,7,9]}}}^{[1,2,3,4,5,6]}   R_{{{10}}}^{[6]}\rfloor\]
\[ R_{{{19}}}^{[9]}= \pi_{{{[1,1,2,3,4,4,5,6,7]}}}^{[1,2,3,4,5,6,7,
8,9]}\lfloor   \pi_{{{[3,4,7]}}}^{[1,2,3]}   R_{{{1}}}^{[3]}, \pi_{{{[
1,2,5,6,8,9]}}}^{[1,2,3,4,5,6]}   R_{{{10}}}^{[6]}\rfloor\]
\[ R_{{{20}}}^{[9]}= \pi_{{{[1,1,2,3,4,4,5,6,7]}}}^{[1,2,3,4,5,6,7,
8,9]}\lfloor   \pi_{{{[3,4,7]}}}^{[1,2,3]}   R_{{{1}}}^{[3]}, \pi_{{{[
1,2,5,6,8,9]}}}^{[1,2,3,4,5,6]}   R_{{{9}}}^{[6]}\rfloor\]
\[ R_{{{21}}}^{[9]}= \pi_{{{[1,2,2,3,4,4,5,6,7]}}}^{[1,2,3,4,5,6,7,
8,9]}\lfloor   \pi_{{{[1,4,8]}}}^{[1,2,3]}   R_{{{1}}}^{[3]}, \pi_{{{[
2,3,5,6,7,9]}}}^{[1,2,3,4,5,6]}   R_{{{9}}}^{[6]}\rfloor\]
\[ R_{{{22}}}^{[9]}= \pi_{{{[1,2,2,3,4,4,5,6,7]}}}^{[1,2,3,4,5,6,7,
8,9]}\lfloor   \pi_{{{[1,4,9]}}}^{[1,2,3]}   R_{{{1}}}^{[3]}, \pi_{{{[
2,3,5,6,7,8]}}}^{[1,2,3,4,5,6]}   R_{{{9}}}^{[6]}\rfloor\]
\[ R_{{{23}}}^{[9]}= \pi_{{{[1,2,2,3,4,4,5,6,7]}}}^{[1,2,3,4,5,6,7,
8,9]}\lfloor   \pi_{{{[1,4,8]}}}^{[1,2,3]}   R_{{{1}}}^{[3]}, \pi_{{{[
2,3,5,6,7,9]}}}^{[1,2,3,4,5,6]}   R_{{{10}}}^{[6]}\rfloor\]
\[ R_{{{24}}}^{[9]}= \pi_{{{[1,2,2,3,4,4,5,6,7]}}}^{[1,2,3,4,5,6,7,
8,9]}\lfloor   \pi_{{{[1,4,9]}}}^{[1,2,3]}   R_{{{1}}}^{[3]}, \pi_{{{[
2,3,5,6,7,8]}}}^{[1,2,3,4,5,6]}   R_{{{10}}}^{[6]}\rfloor\]
\[ R_{{{25}}}^{[9]}= \pi_{{{[1,2,2,3,4,4,5,6,7]}}}^{[1,2,3,4,5,6,7,
8,9]}\lfloor   \pi_{{{[1,4,7]}}}^{[1,2,3]}   R_{{{1}}}^{[3]}, \pi_{{{[
2,3,5,6,8,9]}}}^{[1,2,3,4,5,6]}   R_{{{9}}}^{[6]}\rfloor\]
\[ R_{{{26}}}^{[9]}= \pi_{{{[1,2,2,3,4,4,5,6,7]}}}^{[1,2,3,4,5,6,7,
8,9]}\lfloor   \pi_{{{[1,4,7]}}}^{[1,2,3]}   R_{{{1}}}^{[3]}, \pi_{{{[
2,3,5,6,8,9]}}}^{[1,2,3,4,5,6]}   R_{{{10}}}^{[6]}\rfloor\]
\[ R_{{{27}}}^{[9]}= \pi_{{{[1,1,2,3,4,4,5,6,7]}}}^{[1,2,3,4,5,6,7,
8,9]}\lfloor   \pi_{{{[3,4,9]}}}^{[1,2,3]}   R_{{{1}}}^{[3]}, \pi_{{{[
1,2,5,6,7,8]}}}^{[1,2,3,4,5,6]}   R_{{{9}}}^{[6]}\rfloor\]
\[ R_{{{28}}}^{[9]}= \pi_{{{[1,2,2,3,4,5,6,7,8]}}}^{[1,2,3,4,5,6,7,
8,9]}\lfloor   \pi_{{{[1,7,8]}}}^{[1,2,3]}   R_{{{1}}}^{[3]}, \pi_{{{[
2,3,4,5,6,9]}}}^{[1,2,3,4,5,6]}   R_{{{9}}}^{[6]}\rfloor\]
\[ R_{{{29}}}^{[9]}= \pi_{{{[1,2,2,3,4,5,6,7,8]}}}^{[1,2,3,4,5,6,7,
8,9]}\lfloor   \pi_{{{[1,7,9]}}}^{[1,2,3]}   R_{{{1}}}^{[3]}, \pi_{{{[
2,3,4,5,6,8]}}}^{[1,2,3,4,5,6]}   R_{{{9}}}^{[6]}\rfloor\]
\[ R_{{{30}}}^{[9]}= \pi_{{{[1,2,3,4,5,6,7,7,8]}}}^{[1,2,3,4,5,6,7,
8,9]}\lfloor   \pi_{{{[1,3,9]}}}^{[1,2,3]}   R_{{{1}}}^{[3]}, \pi_{{{[
2,4,5,6,7,8]}}}^{[1,2,3,4,5,6]}   R_{{{9}}}^{[6]}\rfloor\]
\[ R_{{{31}}}^{[9]}= \pi_{{{[1,2,3,4,5,6,7,8,8]}}}^{[1,2,3,4,5,6,7,
8,9]}\lfloor   \pi_{{{[1,3,7]}}}^{[1,2,3]}   R_{{{1}}}^{[3]}, \pi_{{{[
2,4,5,6,8,9]}}}^{[1,2,3,4,5,6]}   R_{{{9}}}^{[6]}\rfloor\]
\[ R_{{{32}}}^{[9]}= \pi_{{{[1,2,3,4,5,6,7,7,8]}}}^{[1,2,3,4,5,6,7,
8,9]}\lfloor   \pi_{{{[2,3,9]}}}^{[1,2,3]}   R_{{{1}}}^{[3]}, \pi_{{{[
1,4,5,6,7,8]}}}^{[1,2,3,4,5,6]}   R_{{{9}}}^{[6]}\rfloor\]
\[ R_{{{33}}}^{[9]}= \pi_{{{[1,2,3,4,5,6,7,8,8]}}}^{[1,2,3,4,5,6,7,
8,9]}\lfloor   \pi_{{{[2,3,7]}}}^{[1,2,3]}   R_{{{1}}}^{[3]}, \pi_{{{[
1,4,5,6,8,9]}}}^{[1,2,3,4,5,6]}   R_{{{9}}}^{[6]}\rfloor\]
\[ R_{{{34}}}^{[9]}= \pi_{{{[1,1,2,3,4,5,6,7,8]}}}^{[1,2,3,4,5,6,7,
8,9]}\lfloor   \pi_{{{[3,7,9]}}}^{[1,2,3]}   R_{{{1}}}^{[3]}, \pi_{{{[
1,2,4,5,6,8]}}}^{[1,2,3,4,5,6]}   R_{{{9}}}^{[6]}\rfloor\]
\[ R_{{{35}}}^{[9]}= \pi_{{{[1,1,2,3,4,5,6,7,8]}}}^{[1,2,3,4,5,6,7,
8,9]}\lfloor   \pi_{{{[3,7,8]}}}^{[1,2,3]}   R_{{{1}}}^{[3]}, \pi_{{{[
1,2,4,5,6,9]}}}^{[1,2,3,4,5,6]}   R_{{{9}}}^{[6]}\rfloor\]
\[ R_{{{36}}}^{[9]}= \pi_{{{[1,2,3,4,5,6,7]}}}^{[1,2,3,4,5,6,7]}
\lfloor   \pi_{{{[2,3,6]}}}^{[1,2,3]}   R_{{{1}}}^{[3]}, \pi_{{{[1,4,5
,7]}}}^{[1,2,3,4]}   R_{{{2}}}^{[6]}\rfloor\]
\[ R_{{{37}}}^{[9]}= \pi_{{{[1,2,3,4,5,6,7]}}}^{[1,2,3,4,5,6,7]}
\lfloor   \pi_{{{[2,3,7]}}}^{[1,2,3]}   R_{{{1}}}^{[3]}, \pi_{{{[1,4,5
,6]}}}^{[1,2,3,4]}   R_{{{2}}}^{[6]}\rfloor\]
\[ R_{{{38}}}^{[9]}= \pi_{{{[1,2,3,4,5,6,7]}}}^{[1,2,3,4,5,6,7]}
\lfloor   \pi_{{{[1,3,7]}}}^{[1,2,3]}   R_{{{1}}}^{[3]}, \pi_{{{[2,4,5
,6]}}}^{[1,2,3,4]}   R_{{{2}}}^{[6]}\rfloor\]
\[ R_{{{39}}}^{[9]}= \pi_{{{[1,2,3,4,5,6,7]}}}^{[1,2,3,4,5,6,7]}
\lfloor   \pi_{{{[1,3,6]}}}^{[1,2,3]}   R_{{{1}}}^{[3]}, \pi_{{{[2,4,5
,7]}}}^{[1,2,3,4]}   R_{{{2}}}^{[6]}\rfloor\]
\[ R_{{{40}}}^{[9]}= \pi_{{{[1,2,3,4,5,6,7]}}}^{[1,2,3,4,5,6,7]}
\lfloor   \pi_{{{[1,5,6]}}}^{[1,2,3]}   R_{{{1}}}^{[3]}, \pi_{{{[2,3,4
,7]}}}^{[1,2,3,4]}   R_{{{1}}}^{[6]}\rfloor\]
\[ R_{{{41}}}^{[9]}= \pi_{{{[1,2,3,4,5,6,7]}}}^{[1,2,3,4,5,6,7]}
\lfloor   \pi_{{{[1,5,7]}}}^{[1,2,3]}   R_{{{1}}}^{[3]}, \pi_{{{[2,3,4
,6]}}}^{[1,2,3,4]}   R_{{{1}}}^{[6]}\rfloor\]
\[ R_{{{42}}}^{[9]}= \pi_{{{[1,2,3,4,5,6,7]}}}^{[1,2,3,4,5,6,7]}
\lfloor   \pi_{{{[2,5,7]}}}^{[1,2,3]}   R_{{{1}}}^{[3]}, \pi_{{{[1,3,4
,6]}}}^{[1,2,3,4]}   R_{{{1}}}^{[6]}\rfloor\]
\[ R_{{{43}}}^{[9]}= \pi_{{{[1,2,3,4,5,6,7]}}}^{[1,2,3,4,5,6,7]}
\lfloor   \pi_{{{[2,5,6]}}}^{[1,2,3]}   R_{{{1}}}^{[3]}, \pi_{{{[1,3,4
,7]}}}^{[1,2,3,4]}   R_{{{1}}}^{[6]}\rfloor\]
\[ R_{{{44}}}^{[9]}= \pi_{{{[1,2,3,4,5,6,7]}}}^{[1,2,3,4,5,6,7]}
\lfloor   \pi_{{{[2,3,6]}}}^{[1,2,3]}   R_{{{1}}}^{[3]}, \pi_{{{[1,4,5
,7]}}}^{[1,2,3,4]}   R_{{{1}}}^{[6]}\rfloor\]
\[ R_{{{45}}}^{[9]}= \pi_{{{[1,2,3,4,5,6,7]}}}^{[1,2,3,4,5,6,7]}
\lfloor   \pi_{{{[2,3,7]}}}^{[1,2,3]}   R_{{{1}}}^{[3]}, \pi_{{{[1,4,5
,6]}}}^{[1,2,3,4]}   R_{{{1}}}^{[6]}\rfloor\]
\[ R_{{{46}}}^{[9]}= \pi_{{{[1,2,3,4,5,6,7]}}}^{[1,2,3,4,5,6,7]}
\lfloor   \pi_{{{[1,3,6]}}}^{[1,2,3]}   R_{{{1}}}^{[3]}, \pi_{{{[2,4,5
,7]}}}^{[1,2,3,4]}   R_{{{1}}}^{[6]}\rfloor\]
\[ R_{{{47}}}^{[9]}= \pi_{{{[1,2,3,4,5,6,7]}}}^{[1,2,3,4,5,6,7]}
\lfloor   \pi_{{{[1,3,7]}}}^{[1,2,3]}   R_{{{1}}}^{[3]}, \pi_{{{[2,4,5
,6]}}}^{[1,2,3,4]}   R_{{{1}}}^{[6]}\rfloor\]
\[ R_{{{48}}}^{[9]}= \pi_{{{[1,2,3,4,5,6,7]}}}^{[1,2,3,4,5,6,7]}
\lfloor   \pi_{{{[1,5,6]}}}^{[1,2,3]}   R_{{{1}}}^{[3]}, \pi_{{{[2,3,4
,7]}}}^{[1,2,3,4]}   R_{{{2}}}^{[6]}\rfloor\]
\[ R_{{{49}}}^{[9]}= \pi_{{{[1,2,3,4,5,6,7]}}}^{[1,2,3,4,5,6,7]}
\lfloor   \pi_{{{[1,5,7]}}}^{[1,2,3]}   R_{{{1}}}^{[3]}, \pi_{{{[2,3,4
,6]}}}^{[1,2,3,4]}   R_{{{2}}}^{[6]}\rfloor\]
\[ R_{{{50}}}^{[9]}= \pi_{{{[1,2,3,4,5,6,7]}}}^{[1,2,3,4,5,6,7]}
\lfloor   \pi_{{{[2,5,7]}}}^{[1,2,3]}   R_{{{1}}}^{[3]}, \pi_{{{[1,3,4
,6]}}}^{[1,2,3,4]}   R_{{{2}}}^{[6]}\rfloor\]
\[ R_{{{51}}}^{[9]}= \pi_{{{[1,2,3,4,5,6,7]}}}^{[1,2,3,4,5,6,7]}
\lfloor   \pi_{{{[2,5,6]}}}^{[1,2,3]}   R_{{{1}}}^{[3]}, \pi_{{{[1,3,4
,7]}}}^{[1,2,3,4]}   R_{{{2}}}^{[6]}\rfloor\]
\[ R_{{{52}}}^{[9]}= \pi_{{{[1,2,3,4,5,6,7,8,9]}}}^{[1,2,3,4,5,6,7,
8,9]}\lfloor   \pi_{{{[2,3,9]}}}^{[1,2,3]}   R_{{{1}}}^{[3]}, \pi_{{{[
1,4,5,6,7,8]}}}^{[1,2,3,4,5,6]}   R_{{{12}}}^{[6]}\rfloor\]
\[ R_{{{53}}}^{[9]}= \pi_{{{[1,2,3,4,5,6,7,8,9]}}}^{[1,2,3,4,5,6,7,
8,9]}\lfloor   \pi_{{{[2,3,8]}}}^{[1,2,3]}   R_{{{1}}}^{[3]}, \pi_{{{[
1,4,5,6,7,9]}}}^{[1,2,3,4,5,6]}   R_{{{12}}}^{[6]}\rfloor\]
\[ R_{{{54}}}^{[9]}= \pi_{{{[1,2,3,4,5,6,7,8,9]}}}^{[1,2,3,4,5,6,7,
8,9]}\lfloor   \pi_{{{[2,3,9]}}}^{[1,2,3]}   R_{{{1}}}^{[3]}, \pi_{{{[
1,4,5,6,7,8]}}}^{[1,2,3,4,5,6]}   R_{{{11}}}^{[6]}\rfloor\]
\[ R_{{{55}}}^{[9]}= \pi_{{{[1,2,3,4,5,6,7,8,9]}}}^{[1,2,3,4,5,6,7,
8,9]}\lfloor   \pi_{{{[2,3,8]}}}^{[1,2,3]}   R_{{{1}}}^{[3]}, \pi_{{{[
1,4,5,6,7,9]}}}^{[1,2,3,4,5,6]}   R_{{{11}}}^{[6]}\rfloor\]
\[ R_{{{56}}}^{[9]}= \pi_{{{[1,2,3,4,5,6,7,8,9]}}}^{[1,2,3,4,5,6,7,
8,9]}\lfloor   \pi_{{{[2,3,7]}}}^{[1,2,3]}   R_{{{1}}}^{[3]}, \pi_{{{[
1,4,5,6,8,9]}}}^{[1,2,3,4,5,6]}   R_{{{12}}}^{[6]}\rfloor\]
\[ R_{{{57}}}^{[9]}= \pi_{{{[1,2,3,4,5,6,7,8,9]}}}^{[1,2,3,4,5,6,7,
8,9]}\lfloor   \pi_{{{[2,3,7]}}}^{[1,2,3]}   R_{{{1}}}^{[3]}, \pi_{{{[
1,4,5,6,8,9]}}}^{[1,2,3,4,5,6]}   R_{{{11}}}^{[6]}\rfloor\]
\[ R_{{{58}}}^{[9]}= \pi_{{{[1,2,3,4,5,6,7,8,9]}}}^{[1,2,3,4,5,6,7,
8,9]}\lfloor   \pi_{{{[1,3,9]}}}^{[1,2,3]}   R_{{{1}}}^{[3]}, \pi_{{{[
2,4,5,6,7,8]}}}^{[1,2,3,4,5,6]}   R_{{{12}}}^{[6]}\rfloor\]
\[ R_{{{59}}}^{[9]}= \pi_{{{[1,2,3,4,5,6,7,8,9]}}}^{[1,2,3,4,5,6,7,
8,9]}\lfloor   \pi_{{{[1,3,8]}}}^{[1,2,3]}   R_{{{1}}}^{[3]}, \pi_{{{[
2,4,5,6,7,9]}}}^{[1,2,3,4,5,6]}   R_{{{12}}}^{[6]}\rfloor\]
\[ R_{{{60}}}^{[9]}= \pi_{{{[1,2,3,4,5,6,7,8,9]}}}^{[1,2,3,4,5,6,7,
8,9]}\lfloor   \pi_{{{[1,3,9]}}}^{[1,2,3]}   R_{{{1}}}^{[3]}, \pi_{{{[
2,4,5,6,7,8]}}}^{[1,2,3,4,5,6]}   R_{{{11}}}^{[6]}\rfloor\]
\[ R_{{{61}}}^{[9]}= \pi_{{{[1,2,3,4,5,6,7,8,9]}}}^{[1,2,3,4,5,6,7,
8,9]}\lfloor   \pi_{{{[1,3,8]}}}^{[1,2,3]}   R_{{{1}}}^{[3]}, \pi_{{{[
2,4,5,6,7,9]}}}^{[1,2,3,4,5,6]}   R_{{{11}}}^{[6]}\rfloor\]
\[ R_{{{62}}}^{[9]}= \pi_{{{[1,2,3,4,5,6,7,8,9]}}}^{[1,2,3,4,5,6,7,
8,9]}\lfloor   \pi_{{{[1,3,7]}}}^{[1,2,3]}   R_{{{1}}}^{[3]}, \pi_{{{[
2,4,5,6,8,9]}}}^{[1,2,3,4,5,6]}   R_{{{12}}}^{[6]}\rfloor\]
\[ R_{{{63}}}^{[9]}= \pi_{{{[1,2,3,4,5,6,7,8,9]}}}^{[1,2,3,4,5,6,7,
8,9]}\lfloor   \pi_{{{[1,3,7]}}}^{[1,2,3]}   R_{{{1}}}^{[3]}, \pi_{{{[
2,4,5,6,8,9]}}}^{[1,2,3,4,5,6]}   R_{{{11}}}^{[6]}\rfloor\]
\[ R_{{{64}}}^{[9]}= \pi_{{{[1,2,3,4,5,6,7,8,9]}}}^{[1,2,3,4,5,6,7,
8,9]}\lfloor   \pi_{{{[1,7,8]}}}^{[1,2,3]}   R_{{{1}}}^{[3]}, \pi_{{{[
2,3,4,5,6,9]}}}^{[1,2,3,4,5,6]}   R_{{{12}}}^{[6]}\rfloor\]
\[ R_{{{65}}}^{[9]}= \pi_{{{[1,2,3,4,5,6,7,8,9]}}}^{[1,2,3,4,5,6,7,
8,9]}\lfloor   \pi_{{{[2,7,9]}}}^{[1,2,3]}   R_{{{1}}}^{[3]}, \pi_{{{[
1,3,4,5,6,8]}}}^{[1,2,3,4,5,6]}   R_{{{12}}}^{[6]}\rfloor\]
\[ R_{{{66}}}^{[9]}= \pi_{{{[1,2,3,4,5,6,7,8,9]}}}^{[1,2,3,4,5,6,7,
8,9]}\lfloor   \pi_{{{[2,7,8]}}}^{[1,2,3]}   R_{{{1}}}^{[3]}, \pi_{{{[
1,3,4,5,6,9]}}}^{[1,2,3,4,5,6]}   R_{{{12}}}^{[6]}\rfloor\]
\[ R_{{{67}}}^{[9]}= \pi_{{{[1,2,3,4,5,6,7,8,9]}}}^{[1,2,3,4,5,6,7,
8,9]}\lfloor   \pi_{{{[1,7,9]}}}^{[1,2,3]}   R_{{{1}}}^{[3]}, \pi_{{{[
2,3,4,5,6,8]}}}^{[1,2,3,4,5,6]}   R_{{{12}}}^{[6]}\rfloor\]
\[ R_{{{68}}}^{[9]}= \pi_{{{[1,2,3,4,5,6,7,8,9]}}}^{[1,2,3,4,5,6,7,
8,9]}\lfloor   \pi_{{{[1,7,9]}}}^{[1,2,3]}   R_{{{1}}}^{[3]}, \pi_{{{[
2,3,4,5,6,8]}}}^{[1,2,3,4,5,6]}   R_{{{10}}}^{[6]}\rfloor\]
\[ R_{{{69}}}^{[9]}= \pi_{{{[1,2,3,4,5,6,7,8,9]}}}^{[1,2,3,4,5,6,7,
8,9]}\lfloor   \pi_{{{[1,7,8]}}}^{[1,2,3]}   R_{{{1}}}^{[3]}, \pi_{{{[
2,3,4,5,6,9]}}}^{[1,2,3,4,5,6]}   R_{{{10}}}^{[6]}\rfloor\]
\[ R_{{{70}}}^{[9]}= \pi_{{{[1,2,3,4,5,6,7,8,9]}}}^{[1,2,3,4,5,6,7,
8,9]}\lfloor   \pi_{{{[3,7,9]}}}^{[1,2,3]}   R_{{{1}}}^{[3]}, \pi_{{{[
1,2,4,5,6,8]}}}^{[1,2,3,4,5,6]}   R_{{{12}}}^{[6]}\rfloor\]
\[ R_{{{71}}}^{[9]}= \pi_{{{[1,2,3,4,5,6,7,8,9]}}}^{[1,2,3,4,5,6,7,
8,9]}\lfloor   \pi_{{{[3,7,8]}}}^{[1,2,3]}   R_{{{1}}}^{[3]}, \pi_{{{[
1,2,4,5,6,9]}}}^{[1,2,3,4,5,6]}   R_{{{12}}}^{[6]}\rfloor\]
\[ R_{{{72}}}^{[9]}= \pi_{{{[1,2,3,4,5,6,7,8,9]}}}^{[1,2,3,4,5,6,7,
8,9]}\lfloor   \pi_{{{[2,7,8]}}}^{[1,2,3]}   R_{{{1}}}^{[3]}, \pi_{{{[
1,3,4,5,6,9]}}}^{[1,2,3,4,5,6]}   R_{{{10}}}^{[6]}\rfloor\]
\[ R_{{{73}}}^{[9]}= \pi_{{{[1,2,3,4,5,6,7,8,9]}}}^{[1,2,3,4,5,6,7,
8,9]}\lfloor   \pi_{{{[2,7,9]}}}^{[1,2,3]}   R_{{{1}}}^{[3]}, \pi_{{{[
1,3,4,5,6,8]}}}^{[1,2,3,4,5,6]}   R_{{{10}}}^{[6]}\rfloor\]
\[ R_{{{74}}}^{[9]}= \pi_{{{[1,2,3,4,5,6,7,8,9]}}}^{[1,2,3,4,5,6,7,
8,9]}\lfloor   \pi_{{{[3,7,8]}}}^{[1,2,3]}   R_{{{1}}}^{[3]}, \pi_{{{[
1,2,4,5,6,9]}}}^{[1,2,3,4,5,6]}   R_{{{10}}}^{[6]}\rfloor\]
\[ R_{{{75}}}^{[9]}= \pi_{{{[1,2,3,4,5,6,7,8,9]}}}^{[1,2,3,4,5,6,7,
8,9]}\lfloor   \pi_{{{[3,7,9]}}}^{[1,2,3]}   R_{{{1}}}^{[3]}, \pi_{{{[
1,2,4,5,6,8]}}}^{[1,2,3,4,5,6]}   R_{{{10}}}^{[6]}\rfloor\]
\[ R_{{{76}}}^{[9]}= \pi_{{{[1,2,3,4,4,5,6]}}}^{[1,2,3,4,5,6,7]}
\lfloor   \pi_{{{[3,6,7]}}}^{[1,2,3]}   R_{{{2}}}^{[4]}, \pi_{{{[1,2,4
,5]}}}^{[1,2,3,4]}   R_{{{3}}}^{[5]}\rfloor\]
\[ R_{{{77}}}^{[9]}= \pi_{{{[1,2,3,3,4,5,6]}}}^{[1,2,3,4,5,6,7]}
\lfloor   \pi_{{{[1,2,5]}}}^{[1,2,3]}   R_{{{1}}}^{[4]}, \pi_{{{[3,4,6
,7]}}}^{[1,2,3,4]}   R_{{{4}}}^{[5]}\rfloor\]
\[ R_{{{78}}}^{[9]}= \pi_{{{[1,2,3,4,5,6,7]}}}^{[1,2,3,4,5,6,7]}
\lfloor   \pi_{{{[4,6,7]}}}^{[1,2,3]}   R_{{{1}}}^{[4]}, \pi_{{{[1,2,3
,5]}}}^{[1,2,3,4]}   R_{{{3}}}^{[5]}\rfloor\]
\[ R_{{{79}}}^{[9]}= \pi_{{{[1,2,3,4,5,6,7]}}}^{[1,2,3,4,5,6,7]}
\lfloor   \pi_{{{[1,2,4]}}}^{[1,2,3]}   R_{{{2}}}^{[4]}, \pi_{{{[3,5,6
,7]}}}^{[1,2,3,4]}   R_{{{4}}}^{[5]}\rfloor\]
\[ R_{{{80}}}^{[9]}= \pi_{{{[1,2,3,4,5,6,7]}}}^{[1,2,3,4,5,6,7]}
\lfloor   \pi_{{{[1,2,4]}}}^{[1,2,3]}   R_{{{1}}}^{[4]}, \pi_{{{[3,5,6
,7]}}}^{[1,2,3,4]}   R_{{{4}}}^{[5]}\rfloor\]
\[ R_{{{81}}}^{[9]}= \pi_{{{[1,2,3,4,5,6,7]}}}^{[1,2,3,4,5,6,7]}
\lfloor   \pi_{{{[4,6,7]}}}^{[1,2,3]}   R_{{{2}}}^{[4]}, \pi_{{{[1,2,3
,5]}}}^{[1,2,3,4]}   R_{{{3}}}^{[5]}\rfloor\]

(i)
The top canonical generating system for $\mathcal R_{10}$ is
given by
$R_{{{1}}}^{[2]}$
and $R_{{{1}}}^{[6]},\ldots, R_{{{12}}}^{[6]}$
and $R_{{{1}}}^{[10]},\ldots, R_{{{8}}}^{[10]}$; where:
\[ R_{{{1}}}^{[10]}= \pi_{{{[1,1,1,2,2]}}}^{[2,3,4,5,6]}\lfloor  {
}\pi_{{{[5,6,7]}}}^{[1,2,3]}   R_{{{1}}}^{[4]}, \pi_{{{[1,2,3,4]}}}^{[
1,2,3,4]}   R_{{{1}}}^{[6]}\rfloor\]
\[ R_{{{2}}}^{[10]}= \pi_{{{[1,2,3]}}}^{[2,3,4]}\lfloor   \pi_{{{[2]}
}}^{[1]}   R_{{{1}}}^{[1]}, \pi_{{{[1,3,4]}}}^{[1,2,3]}   R_{{{1}
}}^{[9]}\rfloor\]
\[ R_{{{3}}}^{[10]}= \pi_{{{[1,2,3]}}}^{[1,2,3]}\lfloor   \pi_{{{[3]}
}}^{[1]}   R_{{{1}}}^{[1]}, \pi_{{{[1,2,4]}}}^{[1,2,3]}   R_{{{1}
}}^{[9]}\rfloor\]
\[ R_{{{4}}}^{[10]}= \pi_{{{[1,1,2,3,4,5,5]}}}^{[1,2,4,5,6,7,8]}
\lfloor   \pi_{{{[5]}}}^{[1]}   R_{{{1}}}^{[1]}, \pi_{{{[1,2,3,4,6,7,8
]}}}^{[1,2,3,4,5,6,7]}   R_{{{4}}}^{[9]}\rfloor\]
\[ R_{{{5}}}^{[10]}= \pi_{{{[1,1,2,3,4,5,5]}}}^{[1,2,3,4,5,6,7]}
\lfloor   \pi_{{{[4]}}}^{[1]}   R_{{{1}}}^{[1]}, \pi_{{{[1,2,3,5,6,7,8
]}}}^{[1,2,3,4,5,6,7]}   R_{{{16}}}^{[9]}\rfloor\]
\[ R_{{{6}}}^{[10]}= \pi_{{{[1,1,2,3,4,4,5,6,6]}}}^{[1,2,3,4,5,6,7,
8,9]}\lfloor   \pi_{{{[4]}}}^{[1]}   R_{{{1}}}^{[1]}, \pi_{{{[1,2,3,5,
6,7,8,9,10]}}}^{[1,2,3,4,5,6,7,8,9]}   R_{{{55}}}^{[9]}\rfloor\]
\[ R_{{{7}}}^{[10]}= \pi_{{{[1,1,2,3,3,4,5,6,6]}}}^{[1,2,4,5,6,7,8,
9,10]}\lfloor   \pi_{{{[7]}}}^{[1]}   R_{{{1}}}^{[1]}, \pi_{{{[1,2,3,4
,5,6,8,9,10]}}}^{[1,2,3,4,5,6,7,8,9]}   R_{{{64}}}^{[9]}\rfloor\]
\[ R_{{{8}}}^{[10]}= \pi_{{{[1,2,3,4,5,6,7,8]}}}^{[1,2,3,4,5,6,7,8]
}\lfloor   \pi_{{{[1,2,3,5]}}}^{[1,2,3,4]}   R_{{{3}}}^{[5]}, \pi_{{{[
4,6,7,8]}}}^{[1,2,3,4]}   R_{{{4}}}^{[5]}\rfloor\]
Remark: Here $R_{{{13}}}^{[6]},\ldots, R_{{{15}}}^{[6]}$ are generated out by the system above.

(j)
The top canonical generating system for $\mathcal R_{11}$ is
given by
$R_{{{1}}}^{[3]}$
and $R_{{{1}}}^{[7]},\ldots, R_{{{18}}}^{[7]}$
and $R_{{{1}}}^{[11]},\ldots, R_{{{67}}}^{[11]}$; where:
\[ R_{{{1}}}^{[11]}= \pi_{{{[1,2,3,4,5]}}}^{[1,2,3,4,5]}\lfloor  {
}\pi_{{{[1,5]}}}^{[1,2]}   R_{{{1}}}^{[2]}, \pi_{{{[2,3,4]}}}^{[1,2,3]
}   R_{{{1}}}^{[9]}\rfloor\]
\[ R_{{{2}}}^{[11]}= \pi_{{{[1,2,3,4,5,6]}}}^{[1,3,4,5,6,7]}
\lfloor \pi_{{{[2,6,7]}}}^{[1,2,3]}   R_{{{1}}}^{[3]}, \pi_{{{[1,3,4,5]}}}^
{[1,2,3,4]}   R_{{{2}}}^{[8]}\rfloor\]
\[ R_{{{3}}}^{[11]}= \pi_{{{[1,2,3,4,5,6]}}}^{[1,2,3,4,5,6]}
\lfloor \pi_{{{[1,2,7]}}}^{[1,2,3]}   R_{{{1}}}^{[3]}, \pi_{{{[3,4,5,6]}}}^
{[1,2,3,4]}   R_{{{3}}}^{[8]}\rfloor\]
\[ R_{{{4}}}^{[11]}= \pi_{{{[1,2,3,4,5,6]}}}^{[1,2,3,4,5,6]}
\lfloor \pi_{{{[1,2,7]}}}^{[1,2,3]}   R_{{{1}}}^{[3]}, \pi_{{{[3,4,5,6]}}}^
{[1,2,3,4]}   R_{{{4}}}^{[8]}\rfloor\]
\[ R_{{{5}}}^{[11]}= \pi_{{{[1,2,3,4,5,6]}}}^{[1,3,4,5,6,7]}
\lfloor \pi_{{{[2,6,7]}}}^{[1,2,3]}   R_{{{1}}}^{[3]}, \pi_{{{[1,3,4,5]}}}^
{[1,2,3,4]}   R_{{{5}}}^{[8]}\rfloor\]
\[ R_{{{6}}}^{[11]}= \pi_{{{[1,2,3,4,5,6]}}}^{[1,2,3,4,5,6]}
\lfloor \pi_{{{[1,2,7]}}}^{[1,2,3]}   R_{{{1}}}^{[3]}, \pi_{{{[3,4,5,6]}}}^
{[1,2,3,4]}   R_{{{2}}}^{[8]}\rfloor\]
\[ R_{{{7}}}^{[11]}= \pi_{{{[1,2,3,4,5,6]}}}^{[2,3,4,5,6,7]}
\lfloor \pi_{{{[1,6,7]}}}^{[1,2,3]}   R_{{{1}}}^{[3]}, \pi_{{{[2,3,4,5]}}}^
{[1,2,3,4]}   R_{{{3}}}^{[8]}\rfloor\]
\[ R_{{{8}}}^{[11]}= \pi_{{{[1,2,3,4,5,6]}}}^{[2,3,4,5,6,7]}
\lfloor \pi_{{{[1,6,7]}}}^{[1,2,3]}   R_{{{1}}}^{[3]}, \pi_{{{[2,3,4,5]}}}^
{[1,2,3,4]}   R_{{{4}}}^{[8]}\rfloor\]
\[ R_{{{9}}}^{[11]}= \pi_{{{[1,2,3,4,5,6]}}}^{[1,2,3,4,5,6]}
\lfloor \pi_{{{[1,2,7]}}}^{[1,2,3]}   R_{{{1}}}^{[3]}, \pi_{{{[3,4,5,6]}}}^
{[1,2,3,4]}   R_{{{5}}}^{[8]}\rfloor\]
\[ R_{{{10}}}^{[11]}= \pi_{{{[1,2,3,3,4,5,6]}}}^{[1,2,3,4,5,6,7]}
\lfloor   \pi_{{{[2,5,6]}}}^{[1,2,3]}   R_{{{1}}}^{[3]}, \pi_{{{[1,3,4
,7]}}}^{[1,2,3,4]}   R_{{{2}}}^{[8]}\rfloor\]
\[ R_{{{11}}}^{[11]}= \pi_{{{[1,2,3,4,4,5,6]}}}^{[1,2,3,4,5,6,7]}
\lfloor   \pi_{{{[2,3,6]}}}^{[1,2,3]}   R_{{{1}}}^{[3]}, \pi_{{{[1,4,5
,7]}}}^{[1,2,3,4]}   R_{{{3}}}^{[8]}\rfloor\]
\[ R_{{{12}}}^{[11]}= \pi_{{{[1,2,3,4,5,6,7]}}}^{[1,2,3,4,5,6,7]}
\lfloor   \pi_{{{[1,4,5]}}}^{[1,2,3]}   R_{{{1}}}^{[3]}, \pi_{{{[2,3,6
,7]}}}^{[1,2,3,4]}   R_{{{4}}}^{[8]}\rfloor\]
\[ R_{{{13}}}^{[11]}= \pi_{{{[1,2,3,4,5,6,7]}}}^{[1,2,3,4,5,6,7]}
\lfloor   \pi_{{{[3,4,7]}}}^{[1,2,3]}   R_{{{1}}}^{[3]}, \pi_{{{[1,2,5
,6]}}}^{[1,2,3,4]}   R_{{{5}}}^{[8]}\rfloor\]
\[ R_{{{14}}}^{[11]}= \pi_{{{[1,2,3,4,5,6,7]}}}^{[1,2,3,4,5,6,7]}
\lfloor   \pi_{{{[1,4,5]}}}^{[1,2,3]}   R_{{{1}}}^{[3]}, \pi_{{{[2,3,6
,7]}}}^{[1,2,3,4]}   R_{{{3}}}^{[8]}\rfloor\]
\[ R_{{{15}}}^{[11]}= \pi_{{{[1,2,3,4,5,6,7]}}}^{[1,2,3,4,5,6,7]}
\lfloor   \pi_{{{[3,4,7]}}}^{[1,2,3]}   R_{{{1}}}^{[3]}, \pi_{{{[1,2,5
,6]}}}^{[1,2,3,4]}   R_{{{2}}}^{[8]}\rfloor\]
\[ R_{{{16}}}^{[11]}= \pi_{{{[1,2,3,4,5,6,7]}}}^{[1,2,3,4,5,6,7]}
\lfloor   \pi_{{{[1,3,4]}}}^{[1,2,3]}   R_{{{1}}}^{[3]}, \pi_{{{[2,5,6
,7]}}}^{[1,2,3,4]}   R_{{{3}}}^{[8]}\rfloor\]
\[ R_{{{17}}}^{[11]}= \pi_{{{[1,2,3,4,5,6,7]}}}^{[1,2,3,4,5,6,7]}
\lfloor   \pi_{{{[3,5,7]}}}^{[1,2,3]}   R_{{{1}}}^{[3]}, \pi_{{{[1,2,4
,6]}}}^{[1,2,3,4]}   R_{{{2}}}^{[8]}\rfloor\]
\[ R_{{{18}}}^{[11]}= \pi_{{{[1,2,3,4,5,6,7]}}}^{[1,2,3,4,5,6,7]}
\lfloor   \pi_{{{[4,5,7]}}}^{[1,2,3]}   R_{{{1}}}^{[3]}, \pi_{{{[1,2,3
,6]}}}^{[1,2,3,4]}   R_{{{2}}}^{[8]}\rfloor\]
\[ R_{{{19}}}^{[11]}= \pi_{{{[1,2,3,4,5,6,7]}}}^{[1,2,3,4,5,6,7]}
\lfloor   \pi_{{{[1,3,5]}}}^{[1,2,3]}   R_{{{1}}}^{[3]}, \pi_{{{[2,4,6
,7]}}}^{[1,2,3,4]}   R_{{{3}}}^{[8]}\rfloor\]
\[ R_{{{20}}}^{[11]}= \pi_{{{[1,2,3,4,5,6,7]}}}^{[1,2,3,4,5,6,7]}
\lfloor   \pi_{{{[4,5,7]}}}^{[1,2,3]}   R_{{{1}}}^{[3]}, \pi_{{{[1,2,3
,6]}}}^{[1,2,3,4]}   R_{{{5}}}^{[8]}\rfloor\]
\[ R_{{{21}}}^{[11]}= \pi_{{{[1,2,3,4,5,6,7]}}}^{[1,2,3,4,5,6,7]}
\lfloor   \pi_{{{[3,5,7]}}}^{[1,2,3]}   R_{{{1}}}^{[3]}, \pi_{{{[1,2,4
,6]}}}^{[1,2,3,4]}   R_{{{5}}}^{[8]}\rfloor\]
\[ R_{{{22}}}^{[11]}= \pi_{{{[1,2,3,4,5,6,7]}}}^{[1,2,3,4,5,6,7]}
\lfloor   \pi_{{{[1,3,4]}}}^{[1,2,3]}   R_{{{1}}}^{[3]}, \pi_{{{[2,5,6
,7]}}}^{[1,2,3,4]}   R_{{{4}}}^{[8]}\rfloor\]
\[ R_{{{23}}}^{[11]}= \pi_{{{[1,2,3,4,5,6,7]}}}^{[1,2,3,4,5,6,7]}
\lfloor   \pi_{{{[1,3,5]}}}^{[1,2,3]}   R_{{{1}}}^{[3]}, \pi_{{{[2,4,6
,7]}}}^{[1,2,3,4]}   R_{{{4}}}^{[8]}\rfloor\]
\[ R_{{{24}}}^{[11]}= \pi_{{{[1,2,3,4,5,6,7]}}}^{[1,2,3,4,5,6,7]}
\lfloor   \pi_{{{[2,5,7]}}}^{[1,2,3]}   R_{{{1}}}^{[3]}, \pi_{{{[1,3,4
,6]}}}^{[1,2,3,4]}   R_{{{2}}}^{[8]}\rfloor\]
\[ R_{{{25}}}^{[11]}= \pi_{{{[1,2,3,4,5,6,7]}}}^{[1,2,3,4,5,6,7]}
\lfloor   \pi_{{{[1,5,7]}}}^{[1,2,3]}   R_{{{1}}}^{[3]}, \pi_{{{[2,3,4
,6]}}}^{[1,2,3,4]}   R_{{{2}}}^{[8]}\rfloor\]
\[ R_{{{26}}}^{[11]}= \pi_{{{[1,2,3,4,5,6,7]}}}^{[1,2,3,4,5,6,7]}
\lfloor   \pi_{{{[1,3,7]}}}^{[1,2,3]}   R_{{{1}}}^{[3]}, \pi_{{{[2,4,5
,6]}}}^{[1,2,3,4]}   R_{{{3}}}^{[8]}\rfloor\]
\[ R_{{{27}}}^{[11]}= \pi_{{{[1,2,3,4,5,6,7]}}}^{[1,2,3,4,5,6,7]}
\lfloor   \pi_{{{[1,3,6]}}}^{[1,2,3]}   R_{{{1}}}^{[3]}, \pi_{{{[2,4,5
,7]}}}^{[1,2,3,4]}   R_{{{3}}}^{[8]}\rfloor\]
\[ R_{{{28}}}^{[11]}= \pi_{{{[1,2,3,4,5,6,7]}}}^{[1,2,3,4,5,6,7]}
\lfloor   \pi_{{{[2,5,7]}}}^{[1,2,3]}   R_{{{1}}}^{[3]}, \pi_{{{[1,3,4
,6]}}}^{[1,2,3,4]}   R_{{{5}}}^{[8]}\rfloor\]
\[ R_{{{29}}}^{[11]}= \pi_{{{[1,2,3,4,5,6,7]}}}^{[1,2,3,4,5,6,7]}
\lfloor   \pi_{{{[1,5,7]}}}^{[1,2,3]}   R_{{{1}}}^{[3]}, \pi_{{{[2,3,4
,6]}}}^{[1,2,3,4]}   R_{{{5}}}^{[8]}\rfloor\]
\[ R_{{{30}}}^{[11]}= \pi_{{{[1,2,3,4,5,6,7]}}}^{[1,2,3,4,5,6,7]}
\lfloor   \pi_{{{[1,3,6]}}}^{[1,2,3]}   R_{{{1}}}^{[3]}, \pi_{{{[2,4,5
,7]}}}^{[1,2,3,4]}   R_{{{4}}}^{[8]}\rfloor\]
\[ R_{{{31}}}^{[11]}= \pi_{{{[1,2,3,4,5,6,7]}}}^{[1,2,3,4,5,6,7]}
\lfloor   \pi_{{{[1,3,7]}}}^{[1,2,3]}   R_{{{1}}}^{[3]}, \pi_{{{[2,4,5
,6]}}}^{[1,2,3,4]}   R_{{{4}}}^{[8]}\rfloor\]
\[ R_{{{32}}}^{[11]}= \pi_{{{[1,2,3,4,4,5,6]}}}^{[1,2,3,4,5,6,7]}
\lfloor   \pi_{{{[1,3,6]}}}^{[1,2,3]}   R_{{{1}}}^{[3]}, \pi_{{{[2,4,5
,7]}}}^{[1,2,3,4]}   R_{{{2}}}^{[8]}\rfloor\]
\[ R_{{{33}}}^{[11]}= \pi_{{{[1,2,3,4,4,5,6]}}}^{[1,2,3,4,5,6,7]}
\lfloor   \pi_{{{[2,3,6]}}}^{[1,2,3]}   R_{{{1}}}^{[3]}, \pi_{{{[1,4,5
,7]}}}^{[1,2,3,4]}   R_{{{2}}}^{[8]}\rfloor\]
\[ R_{{{34}}}^{[11]}= \pi_{{{[1,2,3,3,4,5,6]}}}^{[1,2,3,4,5,6,7]}
\lfloor   \pi_{{{[2,5,7]}}}^{[1,2,3]}   R_{{{1}}}^{[3]}, \pi_{{{[1,3,4
,6]}}}^{[1,2,3,4]}   R_{{{3}}}^{[8]}\rfloor\]
\[ R_{{{35}}}^{[11]}= \pi_{{{[1,2,3,3,4,5,6]}}}^{[1,2,3,4,5,6,7]}
\lfloor   \pi_{{{[2,5,6]}}}^{[1,2,3]}   R_{{{1}}}^{[3]}, \pi_{{{[1,3,4
,7]}}}^{[1,2,3,4]}   R_{{{3}}}^{[8]}\rfloor\]
\[ R_{{{36}}}^{[11]}= \pi_{{{[1,2,3,3,4,5,6]}}}^{[1,2,3,4,5,6,7]}
\lfloor   \pi_{{{[5,6,7]}}}^{[1,2,3]}   R_{{{1}}}^{[4]}, \pi_{{{[1,2,3
,4,8]}}}^{[1,2,3,4,5]}   R_{{{6}}}^{[7]}\rfloor\]
\[ R_{{{37}}}^{[11]}= \pi_{{{[1,2,3,4,5,6]}}}^{[1,2,3,4,5,6]}
\lfloor   \pi_{{{[4,5,6]}}}^{[1,2,3]}   R_{{{1}}}^{[4]}, \pi_{{{[1,2,3,7]}}}
^{[1,2,3,4]}   R_{{{3}}}^{[7]}\rfloor\]
\[ R_{{{38}}}^{[11]}= \pi_{{{[1,2,3,4,5,6]}}}^{[1,3,4,5,6,7]}
\lfloor   \pi_{{{[1,3,4]}}}^{[1,2,3]}   R_{{{2}}}^{[4]}, \pi_{{{[2,5,6,7]}}}
^{[1,2,3,4]}   R_{{{4}}}^{[7]}\rfloor\]
\[ R_{{{39}}}^{[11]}= \pi_{{{[1,2,3,4,4,5,6]}}}^{[1,3,4,5,6,7,8]}
\lfloor   \pi_{{{[1,3,4]}}}^{[1,2,3]}   R_{{{2}}}^{[4]}, \pi_{{{[2,5,6
,7,8]}}}^{[1,2,3,4,5]}   R_{{{5}}}^{[7]}\rfloor\]
\[ R_{{{40}}}^{[11]}= \pi_{{{[1,2,3,4,5,6]}}}^{[1,2,3,4,5,6]}
\lfloor   \pi_{{{[3,5,6]}}}^{[1,2,3]}   R_{{{2}}}^{[4]}, \pi_{{{[1,2,4,7]}}}
^{[1,2,3,4]}   R_{{{3}}}^{[7]}\rfloor\]
\[ R_{{{41}}}^{[11]}= \pi_{{{[1,2,3,4,5,6]}}}^{[2,3,4,5,6,7]}
\lfloor   \pi_{{{[2,3,5]}}}^{[1,2,3]}   R_{{{1}}}^{[4]}, \pi_{{{[1,4,6,7]}}}
^{[1,2,3,4]}   R_{{{4}}}^{[7]}\rfloor\]
\[ R_{{{42}}}^{[11]}= \pi_{{{[1,2,3,3,4,5,6]}}}^{[1,2,3,4,5,6,7]}
\lfloor   \pi_{{{[5,6,7]}}}^{[1,2,3]}   R_{{{2}}}^{[4]}, \pi_{{{[1,2,3
,4,8]}}}^{[1,2,3,4,5]}   R_{{{6}}}^{[7]}\rfloor\]
\[ R_{{{43}}}^{[11]}= \pi_{{{[1,2,3,4,5,6]}}}^{[1,2,3,4,5,6]}
\lfloor   \pi_{{{[4,5,6]}}}^{[1,2,3]}   R_{{{2}}}^{[4]}, \pi_{{{[1,2,3,7]}}}
^{[1,2,3,4]}   R_{{{3}}}^{[7]}\rfloor\]
\[ R_{{{44}}}^{[11]}= \pi_{{{[1,2,3,4,5,6]}}}^{[2,3,4,5,6,7]}
\lfloor   \pi_{{{[2,3,4]}}}^{[1,2,3]}   R_{{{1}}}^{[4]}, \pi_{{{[1,5,6,7]}}}
^{[1,2,3,4]}   R_{{{4}}}^{[7]}\rfloor\]
\[ R_{{{45}}}^{[11]}= \pi_{{{[1,2,3,4,4,5,6]}}}^{[2,3,4,5,6,7,8]}
\lfloor   \pi_{{{[2,3,4]}}}^{[1,2,3]}   R_{{{1}}}^{[4]}, \pi_{{{[1,5,6
,7,8]}}}^{[1,2,3,4,5]}   R_{{{5}}}^{[7]}\rfloor\]
\[ R_{{{46}}}^{[11]}= \pi_{{{[1,2,3,4,5,6,7]}}}^{[1,2,3,4,5,6,7]}
\lfloor   \pi_{{{[1,3,7]}}}^{[1,2,3]}   R_{{{1}}}^{[4]}, \pi_{{{[2,4,5
,6]}}}^{[1,2,3,4]}   R_{{{4}}}^{[7]}\rfloor\]
\[ R_{{{47}}}^{[11]}= \pi_{{{[1,2,3,4,5,6,7]}}}^{[1,2,3,4,5,6,7]}
\lfloor   \pi_{{{[1,3,6]}}}^{[1,2,3]}   R_{{{1}}}^{[4]}, \pi_{{{[2,4,5
,7]}}}^{[1,2,3,4]}   R_{{{4}}}^{[7]}\rfloor\]
\[ R_{{{48}}}^{[11]}= \pi_{{{[1,2,3,4,5,6,7]}}}^{[1,2,3,4,5,6,7]}
\lfloor   \pi_{{{[2,5,7]}}}^{[1,2,3]}   R_{{{2}}}^{[4]}, \pi_{{{[1,3,4
,6]}}}^{[1,2,3,4]}   R_{{{3}}}^{[7]}\rfloor\]
\[ R_{{{49}}}^{[11]}= \pi_{{{[1,2,3,4,5,6,7]}}}^{[1,2,3,4,5,6,7]}
\lfloor   \pi_{{{[1,5,7]}}}^{[1,2,3]}   R_{{{2}}}^{[4]}, \pi_{{{[2,3,4
,6]}}}^{[1,2,3,4]}   R_{{{3}}}^{[7]}\rfloor\]
\[ R_{{{50}}}^{[11]}= \pi_{{{[1,2,3,4,5,6,7]}}}^{[1,2,3,4,5,6,7]}
\lfloor   \pi_{{{[1,2,7]}}}^{[1,2,3]}   R_{{{1}}}^{[4]}, \pi_{{{[3,4,5
,6]}}}^{[1,2,3,4]}   R_{{{4}}}^{[7]}\rfloor\]
\[ R_{{{51}}}^{[11]}= \pi_{{{[1,2,3,4,5,6,7]}}}^{[1,2,3,4,5,6,7]}
\lfloor   \pi_{{{[1,2,6]}}}^{[1,2,3]}   R_{{{1}}}^{[4]}, \pi_{{{[3,4,5
,7]}}}^{[1,2,3,4]}   R_{{{4}}}^{[7]}\rfloor\]
\[ R_{{{52}}}^{[11]}= \pi_{{{[1,2,3,4,5,6,7]}}}^{[1,2,3,4,5,6,7]}
\lfloor   \pi_{{{[1,6,7]}}}^{[1,2,3]}   R_{{{2}}}^{[4]}, \pi_{{{[2,3,4
,5]}}}^{[1,2,3,4]}   R_{{{3}}}^{[7]}\rfloor\]
\[ R_{{{53}}}^{[11]}= \pi_{{{[1,2,3,4,5,6,7]}}}^{[1,2,3,4,5,6,7]}
\lfloor   \pi_{{{[2,6,7]}}}^{[1,2,3]}   R_{{{2}}}^{[4]}, \pi_{{{[1,3,4
,5]}}}^{[1,2,3,4]}   R_{{{3}}}^{[7]}\rfloor\]
\[ R_{{{54}}}^{[11]}= \pi_{{{[1,2,3,4,5,6,7]}}}^{[1,2,3,4,5,6,7]}
\lfloor   \pi_{{{[2,3,7]}}}^{[1,2,3]}   R_{{{1}}}^{[4]}, \pi_{{{[1,4,5
,6]}}}^{[1,2,3,4]}   R_{{{4}}}^{[7]}\rfloor\]
\[ R_{{{55}}}^{[11]}= \pi_{{{[1,2,3,4,5,6,7]}}}^{[1,2,3,4,5,6,7]}
\lfloor   \pi_{{{[2,5,6]}}}^{[1,2,3]}   R_{{{2}}}^{[4]}, \pi_{{{[1,3,4
,7]}}}^{[1,2,3,4]}   R_{{{3}}}^{[7]}\rfloor\]
\[ R_{{{56}}}^{[11]}= \pi_{{{[1,2,3,4,5,6,7]}}}^{[1,2,3,4,5,6,7]}
\lfloor   \pi_{{{[2,3,6]}}}^{[1,2,3]}   R_{{{1}}}^{[4]}, \pi_{{{[1,4,5
,7]}}}^{[1,2,3,4]}   R_{{{4}}}^{[7]}\rfloor\]
\[ R_{{{57}}}^{[11]}= \pi_{{{[1,2,3,4,5,6,7]}}}^{[1,2,3,4,5,6,7]}
\lfloor   \pi_{{{[1,5,6]}}}^{[1,2,3]}   R_{{{2}}}^{[4]}, \pi_{{{[2,3,4
,7]}}}^{[1,2,3,4]}   R_{{{3}}}^{[7]}\rfloor\]
\[ R_{{{58}}}^{[11]}= \pi_{{{[1,2,3,4,5,6,7]}}}^{[2,3,4,5,6,7,8]}
\lfloor   \pi_{{{[2,3,5]}}}^{[1,2,3]}   R_{{{1}}}^{[4]}, \pi_{{{[1,4,6
,7,8]}}}^{[1,2,3,4,5]}   R_{{{5}}}^{[7]}\rfloor\]
\[ R_{{{59}}}^{[11]}= \pi_{{{[1,2,3,4,5,6,7]}}}^{[1,2,3,4,5,6,7]}
\lfloor   \pi_{{{[4,6,7]}}}^{[1,2,3]}   R_{{{2}}}^{[4]}, \pi_{{{[1,2,3
,5,8]}}}^{[1,2,3,4,5]}   R_{{{6}}}^{[7]}\rfloor\]
\[ R_{{{60}}}^{[11]}= \pi_{{{[1,2,3,4,5,6,7,8]}}}^{[1,2,3,4,5,6,7,8
]}\lfloor   \pi_{{{[2,3,8]}}}^{[1,2,3]}   R_{{{1}}}^{[4]}, \pi_{{{[1,4
,5,6,7]}}}^{[1,2,3,4,5]}   R_{{{5}}}^{[7]}\rfloor\]
\[ R_{{{61}}}^{[11]}= \pi_{{{[1,2,3,4,5,6,7,8]}}}^{[1,2,3,4,5,6,7,8
]}\lfloor   \pi_{{{[2,3,7]}}}^{[1,2,3]}   R_{{{1}}}^{[4]}, \pi_{{{[1,4
,5,6,8]}}}^{[1,2,3,4,5]}   R_{{{5}}}^{[7]}\rfloor\]
\[ R_{{{62}}}^{[11]}= \pi_{{{[1,2,3,4,5,6,7,8]}}}^{[1,2,3,4,5,6,7,8
]}\lfloor   \pi_{{{[1,6,7]}}}^{[1,2,3]}   R_{{{2}}}^{[4]}, \pi_{{{[2,3
,4,5,8]}}}^{[1,2,3,4,5]}   R_{{{6}}}^{[7]}\rfloor\]
\[ R_{{{63}}}^{[11]}= \pi_{{{[1,2,3,4,5,6,7,8]}}}^{[1,2,3,4,5,6,7,8
]}\lfloor   \pi_{{{[2,6,7]}}}^{[1,2,3]}   R_{{{2}}}^{[4]}, \pi_{{{[1,3
,4,5,8]}}}^{[1,2,3,4,5]}   R_{{{6}}}^{[7]}\rfloor\]
\[ R_{{{64}}}^{[11]}= \pi_{{{[1,2,3,3,4,4,5,5,6]}}}^{[1,2,3,4,5,6,7
,8,9]}\lfloor   \pi_{{{[1,2,5,6]}}}^{[1,2,3,4]}   R_{{{3}}}^{[5]},{
}\pi_{{{[3,4,7,8,9,10]}}}^{[1,2,3,4,5,6]}   R_{{{9}}}^{[6]}\rfloor\]
\[ R_{{{65}}}^{[11]}= \pi_{{{[1,2,2,3,3,4,4,5,6]}}}^{[1,3,4,5,6,7,8
,9,10]}\lfloor   \pi_{{{[5,6,9,10]}}}^{[1,2,3,4]}   R_{{{4}}}^{[5]},{
}\pi_{{{[1,2,3,4,7,8]}}}^{[1,2,3,4,5,6]}   R_{{{11}}}^{[6]}\rfloor\]
\[ R_{{{66}}}^{[11]}= \pi_{{{[1,2,3,4,4,5,6,6,7]}}}^{[1,2,3,4,5,6,7
,8,9]}\lfloor   \pi_{{{[1,2,3,6]}}}^{[1,2,3,4]}   R_{{{3}}}^{[5]},{
}\pi_{{{[4,5,7,8,9,10]}}}^{[1,2,3,4,5,6]}   R_{{{9}}}^{[6]}\rfloor\]
\[ R_{{{67}}}^{[11]}= \pi_{{{[1,2,2,3,4,4,5,6,7]}}}^{[1,3,4,5,6,7,8
,9,10]}\lfloor   \pi_{{{[5,8,9,10]}}}^{[1,2,3,4]}   R_{{{4}}}^{[5]},{
}\pi_{{{[1,2,3,4,6,7]}}}^{[1,2,3,4,5,6]}   R_{{{11}}}^{[6]}\rfloor\]
Remark: Here
$R_{{{19}}}^{[7]},\ldots, R_{{{28}}}^{[7]}$ are generated out by the system above.

(k)
The top canonical generating system for $\mathcal R_{12}$ is
given by
$R_{{{1}}}^{[4]},R_{{{2}}}^{[4]} $
and $R_{{{1}}}^{[8]},\ldots, R_{{{5}}}^{[8]}$
and $R_{{{1}}}^{[12]},\ldots, R_{{{8}}}^{[12]}$; where:
\[ R_{{{1}}}^{[12]}= \pi_{{{[1,2,3,3,4,5,6,6]}}}^{[1,2,3,4,5,6,7,8]
}\lfloor   \pi_{{{[5]}}}^{[1]}   R_{{{1}}}^{[1]}, \pi_{{{[1,2,3,4,6,7,
8]}}}^{[1,2,3,4,5,6,7]}   R_{{{20}}}^{[11]}\rfloor\]
\[ R_{{{2}}}^{[12]}= \pi_{{{[1,1,2,3,4,4,5,6]}}}^{[1,2,3,4,5,6,7,8]
}\lfloor   \pi_{{{[4]}}}^{[1]}   R_{{{1}}}^{[1]}, \pi_{{{[1,2,3,5,6,7,
8]}}}^{[1,2,3,4,5,6,7]}   R_{{{22}}}^{[11]}\rfloor\]
\[ R_{{{3}}}^{[12]}= \pi_{{{[1,1,2,3,4,5,6,6]}}}^{[1,2,3,4,5,6,7,8]
}\lfloor   \pi_{{{[5]}}}^{[1]}   R_{{{1}}}^{[1]}, \pi_{{{[1,2,3,4,6,7,
8]}}}^{[1,2,3,4,5,6,7]}   R_{{{24}}}^{[11]}\rfloor\]
\[ R_{{{4}}}^{[12]}= \pi_{{{[1,1,2,3,4,5,6,6]}}}^{[1,2,3,4,5,6,7,8]
}\lfloor   \pi_{{{[4]}}}^{[1]}   R_{{{1}}}^{[1]}, \pi_{{{[1,2,3,5,6,7,
8]}}}^{[1,2,3,4,5,6,7]}   R_{{{26}}}^{[11]}\rfloor\]
\[ R_{{{5}}}^{[12]}= \pi_{{{[1,2,3,4,5,6]}}}^{[1,2,3,4,5,6]}
\lfloor \pi_{{{[1,4,5]}}}^{[1,2,3]}   R_{{{1}}}^{[3]}, \pi_{{{[2,3,6]}}}^{[
1,2,3]}   R_{{{1}}}^{[9]}\rfloor\]
\[ R_{{{6}}}^{[12]}= \pi_{{{[1,2,3,4,5,6]}}}^{[1,2,3,4,5,6]}
\lfloor \pi_{{{[1,4,6]}}}^{[1,2,3]}   R_{{{1}}}^{[3]}, \pi_{{{[2,3,5]}}}^{[
1,2,3]}   R_{{{1}}}^{[9]}\rfloor\]
\[ R_{{{7}}}^{[12]}= \pi_{{{[1,2,3,4,5,6]}}}^{[1,2,3,4,5,6]}
\lfloor \pi_{{{[1,3,6]}}}^{[1,2,3]}   R_{{{1}}}^{[3]}, \pi_{{{[2,4,5]}}}^{[
1,2,3]}   R_{{{1}}}^{[9]}\rfloor\]
\[ R_{{{8}}}^{[12]}= \pi_{{{[1,2,3,4,5,6]}}}^{[1,2,3,4,5,6]}
\lfloor \pi_{{{[1,3,5]}}}^{[1,2,3]}   R_{{{1}}}^{[3]}, \pi_{{{[2,4,6]}}}^{[
1,2,3]}   R_{{{1}}}^{[9]}\rfloor\]
Remark: Here $R_{{{6}}}^{[8]},\ldots, R_{{{31}}}^{[8]}$ are generated out by the system above.

(l)
The top canonical generating system for $\mathcal R_{15}$ is
given by
$R_{{{1}}}^{[3]}$
and $R_{{{1}}}^{[7]},\ldots, R_{{{18}}}^{[7]}$
and $R_{{{1}}}^{[11]},\ldots, R_{{{67}}}^{[11]}$
and $R_{{{1}}}^{[15]},\ldots, R_{{{10}}}^{[15]}$; where:
\[ R_{{{1}}}^{[15]}= \pi_{{{[1,2,3,3,4,5,6,7,8,8]}}}^{[1,2,3,4,5,6,
7,8,9,10]}\lfloor   \pi_{{{[1,2,5]}}}^{[1,2,3]}   R_{{{1}}}^{[4]},{
}\pi_{{{[3,4,6,7,8,9,10]}}}^{[1,2,3,4,5,6,7]}   R_{{{50}}}^{[11]}\rfloor\]
\[ R_{{{2}}}^{[15]}= \pi_{{{[1,1,2,3,4,5,6,6,7,8]}}}^{[1,2,3,4,5,6,
7,8,9,10]}\lfloor   \pi_{{{[6,9,10]}}}^{[1,2,3]}   R_{{{2}}}^{[4]},{
}\pi_{{{[1,2,3,4,5,7,8]}}}^{[1,2,3,4,5,6,7]}   R_{{{52}}}^{[11]}\rfloor\]
\[ R_{{{3}}}^{[15]}= \pi_{{{[1,2,3,4,5,6,7,8,9,9]}}}^{[1,2,3,4,5,6,
7,8,9,10]}\lfloor   \pi_{{{[1,2,4]}}}^{[1,2,3]}   R_{{{1}}}^{[4]},{
}\pi_{{{[3,5,6,7,8,9,10]}}}^{[1,2,3,4,5,6,7]}   R_{{{46}}}^{[11]}\rfloor\]
\[ R_{{{4}}}^{[15]}= \pi_{{{[1,2,3,4,5,6,7,8,9,9]}}}^{[1,2,3,4,5,6,
7,8,9,10]}\lfloor   \pi_{{{[1,2,4]}}}^{[1,2,3]}   R_{{{1}}}^{[4]},{
}\pi_{{{[3,5,6,7,8,9,10]}}}^{[1,2,3,4,5,6,7]}   R_{{{50}}}^{[11]}\rfloor\]
\[ R_{{{5}}}^{[15]}= \pi_{{{[1,2,3,4,5,6,7,8,9,9]}}}^{[1,2,3,4,5,6,
7,8,9,10]}\lfloor   \pi_{{{[1,2,5]}}}^{[1,2,3]}   R_{{{1}}}^{[4]},{
}\pi_{{{[3,4,6,7,8,9,10]}}}^{[1,2,3,4,5,6,7]}   R_{{{46}}}^{[11]}\rfloor\]
\[ R_{{{6}}}^{[15]}= \pi_{{{[1,2,3,4,5,6,7,8,9,9]}}}^{[1,2,3,4,5,6,
7,8,9,10]}\lfloor   \pi_{{{[1,2,6]}}}^{[1,2,3]}   R_{{{1}}}^{[4]},{
}\pi_{{{[3,4,5,7,8,9,10]}}}^{[1,2,3,4,5,6,7]}   R_{{{50}}}^{[11]}\rfloor\]
\[ R_{{{7}}}^{[15]}= \pi_{{{[1,1,2,3,4,5,6,7,8,9]}}}^{[1,2,3,4,5,6,
7,8,9,10]}\lfloor   \pi_{{{[5,9,10]}}}^{[1,2,3]}   R_{{{2}}}^{[4]},{
}\pi_{{{[1,2,3,4,6,7,8]}}}^{[1,2,3,4,5,6,7]}   R_{{{52}}}^{[11]}\rfloor\]
\[ R_{{{8}}}^{[15]}= \pi_{{{[1,1,2,3,4,5,6,7,8,9]}}}^{[1,2,3,4,5,6,
7,8,9,10]}\lfloor   \pi_{{{[6,9,10]}}}^{[1,2,3]}   R_{{{2}}}^{[4]},{
}\pi_{{{[1,2,3,4,5,7,8]}}}^{[1,2,3,4,5,6,7]}   R_{{{48}}}^{[11]}\rfloor\]
\[ R_{{{9}}}^{[15]}= \pi_{{{[1,1,2,3,4,5,6,7,8,9]}}}^{[1,2,3,4,5,6,
7,8,9,10]}\lfloor   \pi_{{{[7,9,10]}}}^{[1,2,3]}   R_{{{2}}}^{[4]},{
}\pi_{{{[1,2,3,4,5,6,8]}}}^{[1,2,3,4,5,6,7]}   R_{{{52}}}^{[11]}\rfloor\]
\[ R_{{{10}}}^{[15]}= \pi_{{{[1,1,2,3,4,5,6,7,8,9]}}}^{[1,2,3,4,5,6
,7,8,9,10]}\lfloor   \pi_{{{[7,9,10]}}}^{[1,2,3]}   R_{{{2}}}^{[4]},{
}\pi_{{{[1,2,3,4,5,6,8]}}}^{[1,2,3,4,5,6,7]}   R_{{{48}}}^{[11]}\rfloor\]
\end{theorem}
We remark that the presentation of the generators like above is typically not unique.
For example,
\[R_{{{1}}}^{[11]}=\pi _{{{[1,2,3,4,5]}}}^{[1,2,3,4,5]}  \lfloor\pi _{{{[1
,5]}}}^{[1,2]}  R_{{{1}}}^{[2]},\pi_{{{[2,3,4]}}}^{[1,2,3]}
R_{{{1}}}^{[9]}\rfloor
=\pi _{{{[1,1,1,2,3,4,5]}}}^{[1,2,3,4,5,6,7]}  \lfloor
\pi _{{{[4,5,6]}}}^{[1,2,3]}  R_{{{1}}}^{[3]},\pi _{{{[1,2,3,7]}}}^
{[1,2,3,4]}  R_{{{2}}}^{[8]}\rfloor
.\]
The full and reduced canonical generating systems can also be computed explicitly.
For obvious reasons, we write down only the cardinality of the generating systems here.
\begin{theorem}
\plabel{th:card}
Regarding the cardinalities of the canonical generating systems:
\[
\begin{array}{c||c|c|c|c|c|c|c|c|c|c|c|c}
n&2&3&4&5&6&7&8&9&10&11&12&15\\
\hline
\card\left.\left( \ext^+ \mathcal R_n\right)\right|_{\reg} &2&4&6&42&91&139&211&937&107&945&70&1253\\
\card\left.\left(\ext^{+0} T\mathcal R_n\right)\right|_{\reg}&1&3&3&15&45&59&145&411&65&275&48&387\\
\card\topp_\pi\left.\left( \ext^+ \mathcal R_n\right)\right|_{\reg}&1&1&2&8&16&29&33&89&21&86&15&96\\
\end{array}
\]
\end{theorem}

\begin{remark}
\plabel{rem:gen}
If one interested only in the maximal asymptotical growth, then it is much easier,
as the structure of $\mathcal R_{\mathrm{odd}}=\conv(\mathcal R_9\cup\mathcal R_{15})$
and  $\mathcal R_{\mathrm{even}}=\conv(\mathcal R_{10}\cup\mathcal R_{12})$ is much simpler.
Indeed, $\mathcal R_{\mathrm{odd}}$ is already projectively generated by
$R^{[3]}_{1}$, $R^{[5]}_{1},\ldots,R^{[5]}_{8}$, $R^{[7]}_{7},\ldots,R^{[7]}_{18}$, $R^{[9]}_{2},R^{[9]}_{3}$, $R^{[9]}_{76},\ldots,R^{[9]}_{81}$ (29 elements);
and $\mathcal R_{\mathrm{even}}$ is already projectively generated by
$R^{[2]}_{1}$,
$R^{[4]}_{1},R^{[4]}_{2}$,
$R^{[6]}_{1},R^{[6]}_{2}$,
$R^{[6]}_{5},R^{[6]}_{6}$,
$R^{[6]}_{9},\ldots,R^{[6]}_{12}$,
$R^{[8]}_{1},R^{[10]}_{8}$ (13 elements).
\end{remark}
Then we can discuss
\begin{proof}[(Sketch of) Proofs.]
Let us consider the generating systems $\mathcal G_n$ given in Theorem \ref{th:gen}.
Let $\mathcal R^{\cand}_n=\pgen \mathcal G_n$.
Firstly, we prove that $\mathcal G_n$ is the top canonical generating system of $\mathcal R^{\cand}_n$.
This can be done via computing the other canonical generating systems.
In short, we prove Theorems \ref{th:gen} and \ref{th:card} with respect to  $\mathcal R^{\cand}_n$ instead of $\mathcal R_n$;
except we keep the canonical reduced generating systems $\mathcal G_n^{\mathrm{cr}}$ for further use.
All this can be realized by using standard techniques of linear programming.
Next, we prove the strict containment relations indicated in Theorem \ref{th:Rn} but for $\mathcal R^{\cand}_n$.
We extend the notations $\mathcal G_n, \mathcal R^{\cand}_n,\mathcal G_n^{\mathrm{cr}} $,  for $n\in\mathbb N^\star\setminus\{1,\ldots,12,15\}$
in the way expected from Theorem \ref{th:Rn}.
Due to the construction of the generating elements and \eqref{eq:gen26}, we
immediately see that $\mathcal R^{\cand}_n\subset \mathcal R_n$.
Then we have to prove that $\mathcal R^{\cand}_n= \mathcal R_n$.
For this, it is sufficient to prove that if $A\in\mathcal G_n$ and $B\in\mathcal G_m$,
$n,m\in\{1,\ldots,12,15\}$,
and $a=\card\supp A$,  $b=\card\supp B$,
and $\mathcal I=\{i_1,\ldots, i_a\}$,  $i_1<\ldots< i_a$, $\mathcal J=\{j_1,\ldots, j_b\}$
$j_1<\ldots< j_n$, and $\mathcal I\,\dot\cup\,\mathcal J=\{1,\ldots,a+b\}$, then
\begin{equation}
\lfloor\pi^{[1,\ldots,a]}_{[i_1,\ldots,i_a]}A,\lfloor\pi^{[1,\ldots,b]}_{[j_1,\ldots,j_b]}B  \rfloor
\plabel{eq:surgen}
\end{equation}
is in $\rgen \mathcal G_{a+b}^{\mathrm{cr}}$.
I.~e.~we have to prove that \eqref{eq:surgen} is majorized by a convex combination of relabelings of $\mathcal G_{a+b}^{\mathrm{cr}}$.
This yields a finite number of problems in linear programming, but, due to the large number of cases, this is the
computationally most intensive part of the proof.
By this, Theorems \ref{th:Rn}, \ref{th:gen}, and \ref{th:card} are proven.
Then Theorems \ref{th:rn} and \ref{th:rpn} are sufficient to check through the elements generating sets $\mathcal G_n$,
which is straightforward.
\end{proof}
\begin{remark}
Finding extremal values and coefficients for convex combinations can be done by linear programming.
An accessible exact solver is \verb"qsopt_ex" of   Applegate,  Cook, Dash, Espinoza \cite{ACDE}.
However, here we need linear programming only for qualitative purposes.
Hence corrected solutions obtained from generic numerical solutions are sufficient (at least, in our case).
Ultimately, the computer algebraic software Maple 2015 can handle all computations on a personal computer in a typical configuration.
\end{remark}
\snewpage
\appendix

\section{The number of alternating permutations}
\plabel{sec:asy}
Let $\tilde E_n$ be the number of those permutations $\sigma\in\Sigma_n$ which are alternating.
Here we also count the case $n=0$.
We may consider the exponential generating function
\begin{commentx}
(cf. Flajolet, Sedgewick \cite{FS})
\end{commentx}
\[\tilde E(x)=\sum_{n=0}^\infty \frac{\tilde E_k}{k!}x^k.\]
A more practical  variant is
\[E(x)=\sum_{n=0}^\infty \frac{E_k}{k!}x^k,\]
where $E_0=E_1=1$, and $E_n=\frac12\tilde E_n$, for $n\geq2$.
This yields the exponential generating function of the ``forward'' alternating permutations.
(The coefficients are the [unsigned non-vanishing] Euler numbers.)
\begin{commentx}
(Cf. Stanley \cite{Sta}.)
\end{commentx}
It is well-known (see Andr\'e \cite{A1}/\cite{A2}, cf. \cite{FS}, \cite{Sta}) that
\[E(x)=\tan \left( \frac x2+\frac\pi 4 \right)=\sec x+\tan x.\]
This makes
\[\tilde E(x)=-1-x+2\,\tan \left( \frac x2+\frac\pi 4 \right)=(-1+2\sec x)+(-x+2\tan x).\]

It is easy to see that asymptotically
\[\frac{\tilde E_n}{n!}\sim 4\left(\frac2\pi\right)^{n+1}\]
as $n\rightarrow\infty$.
\begin{commentx}
Indeed,
\[(\tan x)_{2j-1=k}\sim2\left(2/\pi\right)^{k+1};\]
\[(\sech x)_{2j=k}\sim2\left(2/\pi\right)^{k+1};\]
and terms will alternate parity vanish in both cases.
\end{commentx}
Then the expected number of alternating monomials in commutator monomials of
$X_1,\ldots,X_n$ (each variable with multiplicity $1$) averaged over the permutations of the variables is
\[2^{n-1}\frac{\tilde E_n}{n!}\sim \left(\frac4\pi\right)^{n+1}\]
as $n\geq1$, $n\rightarrow\infty$.

\snewpage

\end{document}